\def\year{2022}\relax
\documentclass[letterpaper]{article} % DO NOT CHANGE THIS
\usepackage{aaai22}  % DO NOT CHANGE THIS
\usepackage{times}  % DO NOT CHANGE THIS
\usepackage{helvet}  % DO NOT CHANGE THIS
\usepackage{courier}  % DO NOT CHANGE THIS
\usepackage[hyphens]{url}  % DO NOT CHANGE THIS
\usepackage{graphicx} % DO NOT CHANGE THIS
\usepackage{natbib}  % DO NOT CHANGE THIS AND DO NOT ADD ANY OPTIONS TO IT
\usepackage{caption} % DO NOT CHANGE THIS AND DO NOT ADD ANY OPTIONS TO IT
\usepackage{algorithm}
\usepackage{algorithmic}
\usepackage{threeparttable}
\usepackage{booktabs}

\usepackage{newfloat}
\usepackage{listings}
\DeclareCaptionStyle{ruled}{labelfont=normalfont,labelsep=colon,strut=off} % DO NOT CHANGE THIS
\floatstyle{ruled}
\newfloat{listing}{tb}{lst}{}
\floatname{listing}{Listing}

\title{Convergence Analysis of Randomized SGDA under NC-P{\L} Condition for Stochastic Minimax Optimization Problems}
\author {
    Zehua Liu,\textsuperscript{\rm 1}
    Zenan Li, \textsuperscript{\rm 2}
    Xiaoming Yuan, \textsuperscript{\rm 1}
    Yuan Yao \textsuperscript{\rm 2}
}
\affiliations {
    \textsuperscript{\rm 1} The University of Hong Kong\\
    \textsuperscript{\rm 2} Nanjing University\\
    liuzehua@connect.hku.hk, lizn@smail.nju.edu.cn, 
    xmyuan@hku.hk,
    y.yao@nju.edu.cn
}
\usepackage{bibentry}
\usepackage{amsmath, bm} 
\usepackage{amssymb}
\usepackage{mathtools}
\usepackage{amsthm}
\usepackage{mathrsfs}
\usepackage{subfigure}
\usepackage{svg}
\usepackage{appendix}
\usepackage{array}

\newtheorem{definition}{Definition}
\newtheorem{assumption}{Assumption}
\newtheorem{lemma}{Lemma}
\newtheorem{theorem}{Theorem}
\newtheorem{corollary}{Corollary}

\newtheorem{proposition}{Proposition}

\makeatletter
\newenvironment{breakablealgorithm}
{% \begin{breakablealgorithm}
	\begin{center}
		\refstepcounter{algorithm}% New algorithm
		\hrule height.8pt depth0pt \kern2pt% \@fs@pre for \@fs@ruled
		\renewcommand{\caption}[2][\relax]{% Make a new \caption
			{\raggedright\textbf{\ALG@name~\thealgorithm} ##2\par}%
			\ifx\relax##1\relax % #1 is \relax
			\addcontentsline{loa}{algorithm}{\protect\numberline{\thealgorithm}##2}%
			\else % #1 is not \relax
			\addcontentsline{loa}{algorithm}{\protect\numberline{\thealgorithm}##1}%
			\fi
			\kern2pt\hrule\kern2pt
		}
	}{% \end{breakablealgorithm}
		\kern2pt\hrule\relax% \@fs@post for \@fs@ruled
	\end{center}
}
\makeatother

\def\d{\mathcal{D}}
\def\e{\mathbb{E}}

\def\r{\mathbb{R}}

\begin{document}

\maketitle

\begin{abstract}
    
We introduce a new analytic framework to analyze the convergence of the Randomized Stochastic Gradient Descent Ascent (RSGDA) algorithm for stochastic minimax optimization problems. 
Under the so-called NC-P{\L} condition on one of the variables, our analysis improves the state-of-the-art convergence results in the current literature and hence broadens the applicable range of the RSGDA. 
We also introduce a simple yet effective strategy to accelerate RSGDA , and empirically validate its efficiency on both synthetic data and real data. 
    
\end{abstract}

\section{Introduction}

Minimax optimization plays an essential role in various areas, from classic game theory to contemporary machine learning problems such as generative adversarial networks (GANs)~\cite{goodfellow2014generative}, adversarial training~\cite{goodfellow2014explaining}, multi-agent reinforcement learn~\cite{dai2018sbeed, zhang2021multi}, and online learning~\cite{cesa2006prediction}. %, to name but a few. 
% mathematics, and economics \cite{bacsar1998dynamic, neumann1928theorie, von2007theory}. 
In this paper, we consider the following standard stochastic minimax optimization problem:
\begin{equation} \label{problem: stochastic minimax problem}
\begin{aligned}
    \min_{x \in \r^m} \max_{y \in \r^n} F(x, y) := \e_z [ f(x, y; z) ],
\end{aligned}
\end{equation}
where $x$ and $y$ refer to two agents with $x$ intending to maximize the payoff function $F(x,y)$ and $y$ aiming to minimize it. 
We introduce a random vector $z$ obeying the given distribution $\mathcal{D}$ to represent the stochastic approximation to the payoff function. %\yy{if z is not necessary in introduction, we may consider to introduce it in the next sections.}
% where $x, y$ are variables and $z$ is a random vector obeying the given distribution $\mathcal{D}$. 
For ease of notation, the primal function of this problem is denoted by $\phi (x) := \max_y F(x, y)$.

Considering the fact that recent minimax problems often involve a large number of variables, 
first-order methods, including stochastic gradient descent ascent (SGDA), stochastic gradient descent of max-oracle (SGDmax), and epoch stochastic gradient descent ascent (ESGDA), 
have become the canonical algorithms to solve problem~\eqref{problem: stochastic minimax problem}.
However, SGDA, SGDmax, and ESGDA all contain different drawbacks.
The SGDA algorithm, which alternates between one stochastic gradient ascent step in $y$ and one stochastic gradient descent step in $x$, 
%is widely used for solving problem~\eqref{problem: stochastic minimax problem}. 
has been well-studied in recent years~\citep{chen2021proximal, heusel2017gans, lei2020sgd, lin2020gradient,mescheder2017numerics, nagarajan2017gradient}. 
However, most analysis of SGDA particularly relies on strong assumptions in $F$
(e.g., the strong concavity in $y$), 
and SGDA often cannot work well in practical problems that do not admit such ideal assumptions (even for some simple cases such as $F(x, y) = xy$). 
SGDmax~\cite{jin2020local,lin2020gradient,nouiehed_solving_2019, sanjabi2018convergence} is another well-analyzed algorithm for solving problem (\ref{problem: stochastic minimax problem}). 
Compared with SGDA, the theoretical result~\citep{jin2020local} guarantees that the SGDmax can converge under much milder assumptions (e.g., $F$ is Lipschitz and smooth). 
However, SGDmax requires a maximization step in $y$ instead of the stochastic gradient ascent step, which 
is computationally difficult to achieve in practice.

Compared with SGDA and SGDmax, ESGDA~\cite{goodfellow2014generative, sinha2017certifying, sebbouh2021randomized} is more popular due to its superior empirical performance. 
Elaborately, ESGDA takes a fixed number of stochastic gradient ascent steps in $y$ followed by a stochastic gradient descent step in $x$ during each iteration, and the goal of the ascent steps is to find a good approximation of $y^* (x) := \mathop{\arg\max} F(x ,y)$. 
Despite its popularity, ESGDA is extremely difficult to analyze, and hence there are few theoretical analyses beyond the convex-concave setting. 
For example, the latest analysis is from \citet{yan2020optimal} who considered ESGDA under the condition that $F$ is weakly convex in $x$ and strongly concave in $y$. 
%To alleviate the difficulty of 

To better analyze ESGDA, a randomized version of ESGDA called RSGDA is proposed~\citep{sebbouh2021randomized} to bridge the theoretical framework and the empirical result. 
Specifically, at each iteration, 
RSGDA takes a stochastic gradient descent step in $x$ with probability $p$ and a stochastic gradient ascent step in $y$ with probability $1-p$. 
Intuitively, RSGDA is consistent with ESGDA in the sense of expectation: during multiple iterations, 
it takes one gradient descent step, followed by $\frac{1-p}{p}$ gradient ascent steps on average.
% Therefore, one may focus on the theoretical analyses of this randomized version RSGDA in the hope of extending them to the deterministic version ESGDA. %\yy{delete this sentence?}

% By introducing the coin toss parameter $p$, theoretical analysis of RSGDA is simple and transparent.
However, the current analysis of RSGDA in \citep{sebbouh2021randomized} is still unsatisfactory.
First, the analysis only provides a partial convergence result that RSGDA can converge to a stationary point of $\phi$, other than the original function $F$. 
In other words, it proves the convergence of RSGDA with respect to $x$, 
but lacks an analysis of $y$ and, more importantly, the joint variable $(x,y)$. 
Second, the provided convergence rate of RSGDA needs to be re-determined. 
Their theoretical result indicates that RSGDA is slower than SGDA, which is inconsistent with numerical experiments showing that RSGDA is at least as fast as SGDA. 
Third, their analysis is limited to the strongly concave setting, which is prohibitively impractical in most cases.

To this end, we propose a new technical framework to analyze RSGDA. 
Elaborately, inspired by \citet{yang2021faster}, we introduce a new Lyapunov function $V$ to bridge the gap between the original function $F(x,y)$ and the primal function $\phi(x)$.
Furthermore, we use $V$ to analyze RSGDA in a relatively more moderate condition, i.e., the NC-P{\L} setting (nonconvex in $x$ and P{\L} condition in $y$), 
and provide more sound convergence results based on our framework. 
In addition, we also analyze the convergence rate with respect to the parameter $p$, which guides a new selection strategy of $p$. %\yy{can RSGDA with you $p$ selection outperform ESGDA? if yes, I think we can have some much stronger statements here.}

% consider RSGDA under the NC-P{\L} setting, i.e., $F(x, y)$ is nonconvex in $x$ and has the P{\L} condition in $y$. 
% which means that $F(x, y)$ is nonconvex in $x$ and has the P{\L} condition in $y$.
% Hence, our work naturally covers the work \citep{sebbouh2021randomized}.

Our contributions can be summarized as follows:

\begin{itemize}
    % \item In Section \ref{sec: rsgda}, we show the convergence results for RSGDA in both the strongly convex strongly concave case (SCSC) and the NC-P{\L} case. Our contributions are threefold. First, we generalize the nonconvex-strongly concave condition considered in \cite{sebbouh2021randomized} to the NC-P{\L} condition. Moreover, different from the convergence analyses in \cite{jin2020local, sebbouh2021randomized}, we prove that we can obtain a pair $(x, y)$ from the sequence generated by RSGDA such that $ \| \nabla F (x, y)\| \leq \epsilon $ and $ \| \nabla \phi (x) \| \leq \epsilon$ simultaneously. Furthermore, we provide convergence results that are better than the ones provided in \cite{sebbouh2021randomized}.
    
    % Compared to aforementioned previous work, our analysis does not assume uniformly strongly concave in $y$, and improves the convergence bound.
    \item We introduce a new framework for the analysis of RSGDA, and prove the almost sure convergence of RSGDA for both $x$ and the joint variable $(x, y)$ in NC-P{\L} setting.
    
    \item We further analyze the convergence rate of RSGDA, and derive a more sound result than previous analysis. Based on this analysis, we also propose a simple but effective method to adjust the parameter $p$ for RSGDA.
    
    \item Empirical experiments show the efficiency of RSGDA and confirm our theoretical results.
\end{itemize}

The rest of this paper is organized as follows. 
Section~\ref{sec:related} briefly overviews related work in the direction of solving stochastic minimax optimization problem. 
Section \ref{sec: preliminaries} is devoted to preliminaries. 
In Section \ref{sec: rsgda}, we analyze RSGDA under the NC-P{\L} assumption. 
Moreover, we propose an intuitive method to determine the parameter $p$.
Section \ref{sec: experiment} contains the numerical experiments. 

\section{Related work} \label{sec:related}

\textbf{P{\L} condition in minimax optimization}. P{\L} condition, named after Polyak and {\L}ojasiewicz, was initially introduced by Polyak in \citet{polyak1963gradient} to obtain the global convergence of gradient descent at a linear rate. The P{\L} condition, roughly speaking, describes the sharpness of a function up to a representation. A generalized form of this condition, which nowadays is called the Kurdyka-{\L}ojasiewicz condition, was introduced by Kurdyak and {\L}ojasiewicz in \citet{kurdyka1998gradients} and \citet{lojasiewicz1963propriete}. There is tremendous work related to the K{\L} condition, making it impossible to list all of them. Curious readers can refer to \citet{bolte2007clarke, bolte2007lojasiewicz} for a systematical discussion. 
In the deterministic case, \citet{nouiehed_solving_2019} showed that GDA and its multi-step variant can achieve an approximate critical point in $O(\epsilon^{-2})$ steps. 
Recently, \citet{fiez2021global} proved that GDA converges to an approximate differential Stackelberg equilibrium with complexity $O(\epsilon^{-2})$; \citet{yang2021faster} proved that a single loop GDA converges to an approximate Stackelberg equilibrium, and an approximate local Nash equilibrium can be constructed from GDA. 

\textbf{Other minimax optimization}. Minimax optimization problems have received wide attention since the work of von Neumann \cite{neumann1928theorie}. Since then, minimax problems have been well studied in the convex-concave setting. However, results beyond the convex-concave setting are much more recent. A large body of existing work~\cite{fiez2021global, lin2020gradient, boct2020alternating, lin2020near} considered GDA in the nonconvex strongly concave setting, obtaining $O(\epsilon^{-2})$ computation complexity in the deterministic case and $O(\epsilon^{-4})$ complexity in the stochastic case. \citet{daskalakis2018limit} considered GDA in the nonconvex-nonconcave setting and provided theoretical analysis of the limit points of GDA. Furthermore, they proposed Optimistic-GDA to robustify the performance of GDA. Due to the difficulty of finding an approximate Nash equilibrium in the general nonconvex-nonconcave setting, different notions of local optimal solutions as well as their properties have been investigated in \citet{jin2020local, fiez2021global} and so on.
\citet{xian_2021_faster} considered the minimax problems in the nonconvex-strongly concave setting and proposed a decentralized algorithm to solve this problem, which achieves a faster convergence rate than SGDA.
\citet{sharma_2022_federated} systemically discussed the local performance of SGDA in both the convex-concave setting and the nonconvex-nonconcave setting.
In \citet{luo_2020_advances}, the authors proposed a variant of SGDA called SREDA, and proved that it achieved the best known stochastic gradient complexity in the nonconvex-strongly concave setting.
Some other work including \citet{diakonikolas_2021_efficient, li_2022_convergence, lee2021fast} considered the applications of extra-gradient method to the minimiax problems.

\section{Preliminaries} \label{sec: preliminaries} 

% In this section, we summarize some preliminaries, define some notations for simplifying the presentation of further analysis, and give some basic results.

\subsection{Notations}

Throughout this paper, we let $\| \cdot \| := \sqrt{\langle \cdot, \cdot \rangle}$ denote the $\ell_2$ norm, and $\langle \cdot, \cdot \rangle$ denote the inner product in the Euclidean space. 
% For any operator $A :\r^m \to \r^n$, $JA$ is denoted by the Jacobian of $A$.
We are interested in the minimax problem of this form: 
\begin{equation} \label{problem: minimax problem}
    \min_{x \in \r^m} \max_{y \in \r^n} F(x, y) := \e_z [f(x, y; z)],
\end{equation}
where $z$ is a random variable obeying a distribution $\d$ and $F$ is nonconvex in $x$ for any fixed $y$ and possibly nonconcave in $y$.
Following \citep{jin2020local}, we define $\phi (x) := \arg\max_y F(x, y)$.
Function $\phi$ plays a bridge between the inner problem in $y$ and the outer problem in $x$.

% For ease of notation, we introduce some notations. Define 
% \begin{equation*}
%     u = \begin{pmatrix}
%     x \\ 
%     y
%     \end{pmatrix} \quad 
%     G(u; \alpha, \eta) = \begin{pmatrix}
%     \alpha \nabla_x F(x, y) \\
%     - \eta \nabla_y F(x, y)
%     \end{pmatrix}
% \end{equation*}

\begin{definition}[Smooth function] \label{def: smooth}
We say a function $g$ is $L$-smooth with $L \geq 0$, if it is differentiable and its gradient $\nabla f$ is $L$-Lipschitz continuous.
\end{definition}

\begin{definition}[$\mu$-strongly convexity]
A differentiable function $g : \r^d \to \r$ is called $\mu$-strongly convex if  
\begin{equation}
    g (y) \geq g(x) + \langle \nabla g(x), y - x \rangle + \frac{\mu}{2} \| y - x \|^2.
\end{equation}
\end{definition}
One can easily extend this concept to the minimax problem, where $F(x, y)$ is called $\mu$-strongly-convex-strongly-concave (SCSC), if $F(\cdot, y)$ is $\mu$-strongly convex for any fixed $y$ and  $-F(x, \cdot)$ is $\mu$-strongly convex for any fixed $x$. 

% \lizn{I think here we need max-oracle definition and some illustrations for it.}

% \begin{definition}[Max-oracle]
% We define function $\phi (x) := \max_y F(x, y)$. Then problem (\ref{problem: minimax problem}) is equivalent to find the minimization of function $\phi$.
% \end{definition}
% \yy{this definition is kind of incomplete to me. Also, if the definitions are from existing work, you may need to insert the references.}
% \begin{definition}
    
% For a differentiable function $F$, 
% \begin{itemize}
%     \item A point $(x^*, y^*)$ is a critical point of $F$ if $\nabla F(x^*, y^*) = 0$.
%     \item A critical point $(x^*, y^*)$ is a local minimax point if there exists a neighborhood $U$ around $(x^*, y^*)$ so that for all $(x, y) \in U$ we have that $F(x^*, y) \leq F(x^*, y^*) \leq (x, y^*)$.
% \end{itemize}

% \end{definition}

\subsection{Optimality}
Generally, there are two reasonable solutions to problem~\eqref{problem: stochastic minimax problem} worthy to be noted. 
In simultaneous games, 
one often seeks a Nash equilibrium $(x^*, y^*)$, in which $x^*$ is a global minimum of $F(\cdot, y^*)$ and $y^*$ is a global maximum of $F(x^*, \cdot)$. 
On the other hand, many recent machine learning tasks focus on sequential games, e.g., adversarial training and generative adversarial network, and aim to achieve a Stackelberg equilibrium $(x^*, y^*)$, in which $x^*$ is a global minimum of $\phi(\cdot)$ and $y^*$ is a global maximum of $F(x^*, \cdot)$. 
The Stackelberg equilibrium is also called the global minimax point in some literature~\citep{jin2020local}. 

However, most of the minimax problems arising in machine learning applications are nonconvex in $x$ and nonconcave in $y$, making finding either a Nash equilibrium or a Stackelberg equilibrium impractical. 
This motivates the quest to propose two different notions of local optimality, i.e., local Nash  equilibrium~\citep{daskalakis2018limit, mazumdar2020gradient} and local Stackelberg  equilibrium~\citep{jin2020local}. 
In a nutshell, 
a point $(x^*, y^*)$ is called a local Nash equilibrium if $x^*$ is a local minimum of $F(\cdot, y^*)$ and $y^*$ is a local maximum of $F(x^*, \cdot)$. 
Similarly,  a point $(x^*, y^*)$ is called a local Stackelberg equilibrium if $x^*$ is a local minimum of $\phi (\cdot)$ and $y^* \in \arg\max_y F(x^*, y)$.

However, 
in nonconvex and nonconcave setting, %(or large-scale setting?)
verifying such local equilibrium is still extremely hard since it often requires second-order optimality condition. 
Hence, we instead find a solution that satisfies first-order necessary condition: 
a Stackelberg-type stationary point $(x^*, y^*)$ if $\nabla \phi (x^*) = 0$ and $y^* \in \arg\min_y F(x^*, y)$, and a Nash-type stationary point $(x^*, y^*)$ if $\nabla F(x^*, y^*) = 0$. 
It should be noted that these two kinds of stationary points may not coincide. For example, considering $F(x,y)=xy+x^3$ with $y \in [-2,2]$, one can observe that $(0,0)$ is only the Nash-type stationary point but not the Stackelberg-type one.

% A straightforward observation indicates that a local Stackelberg equilibrium is a Stackelberg type stationary point and a local Nash equilibrium is a Nash type stationary point. 

\subsection{Assumptions}

The following two assumptions are effective throughout, which are standard in stochastic optimization.

\begin{assumption}[Smoothness] \label{assum: smoothness}

The objective function $F$ is $L_1$-smooth.

\end{assumption}

\begin{assumption}[Sampling rule] \label{assum: bounded variance}
For any $x, y$, the gradient estimation of $F$ is unbiased: 
\begin{equation}
    \e_z [\nabla f(x,y;z)] = \nabla F(x, y).
\end{equation}
and its variance is bounded, i.e., there exists a positive constant $\sigma$, such that for all $(x, y)$, 
\begin{equation}
    \e_z \left[ \| \nabla F(x, y) - \nabla f(x, y; z) \|^2 \right] \leq \sigma^2.
\end{equation}
\end{assumption}

Besides the assumptions above, we further assume the following property of the variable $y$.

\begin{assumption}[P{\L} condition in $y$] \label{assum: ncpl condtion}
For any fixed $x$, the maximization $\max_{y} F(x, y)$ has a nonempty solution set and a finite optimal value. Furthermore, there exists $\mu > 0$ such that 
\begin{equation}
    \| \nabla_y F(x, y) \|^2 \geq 2 \mu [\max_y F(x, y) - F(x, y)]
\end{equation}
holds for all $y$.
\end{assumption}

% The P{\L} condition was originally introduced in \citet{polyak1963gradient} for the purpose to guarantee global convergence of gradient descent at a linear rate. 
% Roughly speaking, P{\L} condition describes the sharpness of a function $g$ up to a reparametrization~\cite{lojasiewicz1963propriete}.\yy{the above sentence has appeared in related work}
It should be noted P{\L} condition is an independent property of convexity. 
In other words, there exists nonconvex function $g$ satisfies the P{\L} condition. 
However, the P{\L} condition holds for any strongly convex function $g$, 
and thus P{\L} condition can be viewed as a non-trivial generalization of strong convexity.
 
% On the other hand, the P{\L} condition is a special case of the famous {\L}ojasiewicz inequality. 
% The {\L}ojasiewicz inequality, 
% % which is named after S.{\L}ojasiewicz (see \cite{lojasiewicz1963propriete}), 
% roughly speaking, describes the sharpness of a function $f$ up to a reparametrization~\cite{lojasiewicz1963propriete}. 
% A generalized form of this inequality, which nowadays is called Kurdyka-{\L}ojasiewicz property, has been introduced by Kurdyka in \cite{kurdyka1998gradients}. \lizn{Then?} 
% The reader can refer to \cite{bolte2007lojasiewicz, bolte2007clarke} for non-smooth cases. 
% K\L-property has been successfully used in various fields of mathematics: minimization and algorithms \cite{absil2005convergence, attouch2009convergence, attouch2010proximal, attouch2013convergence, haraux1998convergence, huang2001convergence}; asymptotic theory of partial differential equations \cite{chill2003convergence, jendoubi1998simple, simon1983asymptotics}; among others.

A direct result derived by P{\L} condition is the smoothness of $\phi$.  
Generally, although $F$ is a smooth function, $\phi$ is not so, and even may not be differentiable. 
In this case, it is nearly impossible to define the stationary point of $\phi$. 
However, by combining Assumption~\ref{assum: smoothness} and Assumption~\ref{assum: ncpl condtion}, we can obtain that $\phi$ is $L_2$-smooth with $L_2 = L_1 + \frac{L_1 \kappa}{2}$ and $\kappa = L_1 / \mu$. 
Here $\kappa$ is also referred to as the conditional number of problem~\eqref{problem: minimax problem} \citep{nouiehed_solving_2019}.
% Precisely, we have the following result.

% \begin{lemma}[\citet{nouiehed_solving_2019}] \label{lemma: 1}
% Under  Assumption \ref{assum: smoothness} and \ref{assum: ncpl condtion}, $\phi$ is $L_2$ smooth with $L_2 = L_1 + \frac{L_1 \kappa}{2}$ and $\kappa = L_1 / \mu$.
% \end{lemma}

\subsection{SGDA and its variants}

In this subsection, we provide the formal formulations of SGDA, SGDmax, ESGDA, and RSGDA in Algorithm \ref{algo: rsgda}.
% First, we recall the recursive iterations of well-known SGDA, SGDmax, and RSGDA.
For details, 
at the \textit{k}-th iteration, $z_k$ is first sampled from $\mathcal{D}$. 
Next, 
SGDA takes one descent step in $x$ along $\nabla_x f(x_k, y_k; z_k)$ and one ascent step in $y$ along $\nabla_x f(x_k, y_k; z_k)$; 
SGDmax takes the descent step in $x$ along $\nabla_x f(x_k, y_k; z_k)$ and calculates the maximum of $y$ for $f(x_k, y; z_k)$; 
ESGDA takes $m$ gradient ascent steps in $y$ and one gradient descent step in $x$.
For RSGDA, it 
% We now introduce RSGDA in detail. 
% At each iteration $k$, we first sample $z_k \sim \d$. 
% Then, we 
takes the descent step in $x$ along the stochastic gradient $\nabla_x f(x_k, y_k; z_k)$ with probability $p$ or the gradient ascent step in $y$ along the stochastic gradient $\nabla_x f(x_k, y_k; z_k)$ with probability $1 - p$. 
% The details is summarized in Algorithm~\ref{algo: rsgda}.

\begin{breakablealgorithm}
\caption{SGDmax/ SGDA/ ESGDA/ RSGDA}  \label{algo: rsgda}

\textbf{Inputs:} initial points $x_0, y_0$, step sizes $\{ (\alpha_k, \eta_k) \}_{k=1}^\infty$, loop size $m$, max-oracle accuracy $\delta$, constant parameter $p$.

\begin{algorithmic}
\FOR{$k = 0, 1, 2, \dots$}
\STATE \underline{\textbf{SGDmax}}:

\quad Find $y_{k+1}$, s.t. $F(x_k, y_{k+1}) \geq \phi (x_k) + \delta$;

\quad Sample $z_k \sim \mathcal{D}$;

\quad $x_{k+1} = x_k - \alpha_k \nabla_x f(x_k, y_{k+1}; z_k)$;
\STATE \underline{\textbf{SGDA}}:

\quad Sample $z_k \sim \mathcal{D}$;

\quad $y_{k+1} = y_k + \eta_k \nabla_y f(x_k, y_k; z_k)$;

\quad $x_{k+1} = x_k - \alpha_k \nabla_x f(x_k, y_{k+1}; z_k)$;
\STATE \underline{\textbf{ESGDA}}:
\FOR{$t=0, \dots, m$}
\STATE Sample $z_k^t \sim \mathcal{D}$;
\STATE $y_k^{t+1} = y_k + \eta_k \nabla_y f(x_k, y_k^t; z_k^t)$;
\ENDFOR

\quad $y_{k+1} = y_k^{m+1}$

\quad Sample $z_k \sim \mathcal{D}$;

\quad $x_{k+1} = x_k - \alpha_k \nabla_x f(x_k, y_{k+1}; z_k)$;

\STATE \underline{\textbf{RSGDA}}:

\quad Sample $z_k \sim \mathcal{D}$;

\quad With probability $p$: 
\begin{equation*}
    x_{k+1} = x_k - \alpha_k \nabla_x f(x_k, y_k; z_k); 
\end{equation*}

\quad With probability $1-p$:
\begin{equation*}
y_{k+1} = y_k + \eta_k \nabla_y f (x_k, y_k; z_k);
\end{equation*}
% \begin{equation*}
% (x_{k+1}, y_{k+1}) = 
% \begin{dcases}
%     (x_k^+, y_k), & \quad \text{w.p.} \quad p, \\ 
%     (x_k, y_k^+), & \quad \text{w.p.} \quad 1 - p.
% \end{dcases}
% \end{equation*}

\ENDFOR
\end{algorithmic}

\end{breakablealgorithm}

\section{Randomized stochastic gradient descent ascent (RSGDA)} \label{sec: rsgda}
%\yy{A table summarizing the conditions and results in introduction would help readers quickly get a big picture of your work.}
In this section, 
we first discuss the motivation of RSGDA (Section \ref{sec: 4.1}),
and show that RSGDA converges to the unique Nash equilibrium of $F$ in the SCSC setting (Section \ref{sec: 4.2}).
Since it is impossible to discuss the convergence to Nash equilibriums in the nonconvex-nonconcave setting, we next show that RSGDA converges to a Nash-type stationary point  and a Stackelberg-type stationary point in the NC-P{\L} condition (Section \ref{sec: 4.3}).
Finally, based on the theoretical results, we propose a selection strategy for $p$ (Section \ref{sec: 4.4}).
    
\subsection{Motivation of RSGDA} \label{sec: 4.1}

The intuition of ESGDA is straightforward. 
On the one hand, SGDmax contains a complete theory but is impractical in applications due to the calculation of $y^*$. 
On the other hand, although SGDA is tractable, it lacks theoretical guarantees and often fails in many cases.
ESGDA takes advantage of both SGDmax and RSGDA. 
% ESGDA takes a fixed number of steps in $y$ and one step in $x$ in each iteration. 
In iterations of ESGDA, 
the multiple steps in $y$ not only provide a reasonable estimation of $y^*$ but are also tractable.

Though ESGDA performs better than SGDA in applications, its theoretical properties are still unclear. 
Analyzing its properties is challenging due to the following two technical reasons.
First, the multiple gradient steps in $y$ create several immediate variables, causing the gap between $y_k$ and $y_{k+1}$. 
Second, a theoretical analysis of the best inner update steps seems to be unachievable if it is analyzed by classical techniques.
Therefore, we focus on the randomized version of ESGDA, i.e., RSGDA, as a surrogate. 
% \lizn{Need a further conclusion sentence}

\subsection{Convergence under SCSC condition} \label{sec: 4.2} %\label{subsec: 3.1}
The existing theoretical analysis~\cite{farnia2021train} demonstrates that, 
when $F$ is SCSC and with noiseless gradients, SGDA and ESGDA can successfully strongly converge to the unique Nash equilibrium with linear convergence rate. 
Thus, we first show that, RSGDA enjoys the same convergence property with SGDA and ESGDA in the SCSC setting. 
In other words, the randomized update of $x$ and $y$ will not damage the convergence results. 
Throughout this subsection, we always assume that $\sigma = 0$, i.e., the gradient estimation is exact.
% Before further details about Algorithm \ref{algo: rsgda}, we first consider it in a simple case: $F$ is $\mu$-SCSC. As discussed above, if $F$ is $\mu$-SCSC, it must have the NC-P{\L} condition. However, due to the rich theoretical results for the SCSC function, it is meaningful to consider Algorithm \ref{algo: rsgda} under the SCSC condition separately. 
Note that when $F$ is SCSC, there is a unique Stackelberg equilibrium and a unique Nash equilibrium of $F$, and the Stackelberg equilibrium and Nash equilibrium are equal~\citep{rockafellar_convex_1970}.
For simplicity, we eliminate the random term $z$ due to the noiseless gradients setting, and assume that the step sizes are constant, i.e.,  for any $k \geq 0$, $\alpha_k \equiv \alpha$ and $\eta_k \equiv \eta$ for some given values $\alpha$ and $\eta$.  
The convergence result is concluded in the following theorem, with all the proofs deferred to the appendix. 
% For ease of notation, we define the operator $H_{\text{RSGDA}}$ as follows: 
% \begin{equation*}
%     H_{\text{RGDA}} (x, y) := 
%     \begin{dcases}
%         (x - \alpha \nabla_x F(x, y), y), & \quad \text{w.p.} \quad p, \\ 
%         (x, y + \eta \nabla_y F(x, y)), & \quad \text{w.p.} \quad 1 - p.
%     \end{dcases}
% \end{equation*}

\begin{theorem} \label{thm: rgda is contractive}
Assume that $F$ is $\mu$-SCSC and $\sigma = 0$. Let $\{ (x_k, y_k) \}$ be the sequence generated by Algorithm \ref{algo: rsgda}, $(x^*, y^*)$ be the Nash equilibrium of $F$. 
For sufficiently small $\alpha = \eta$, there exists a constant $\rho < 1$, such that $\forall k \geq 0$, 
\begin{equation*}
    \e_k [\| (x_{k+1}, y_{k+1}) - (x^*, y^*) \|^2] \leq \rho \| (x_k, y_k) - (x^*, y^*) \|^2.
\end{equation*}
\end{theorem}
% \lizn{Conclusion: vs. GDA}

\noindent{\em Remarks}.
Intuitively, Theorem \ref{thm: rgda is contractive} states that RSGDA is a quasi-contractive operator in the expectation viewpoint. 
Consequently, we can prove that RSGDA converges linearly to the minimax point for a $\mu$-SCSC function $F$ in expectation, which shares the same convergence rate as SGDA and ESGDA in the SCSC setting with $\sigma = 0$. 

\begin{corollary}
    
Consider the setting of Theorem \ref{thm: rgda is contractive}, $\{(x_k, y_k)\}$ converges to the minimax point $(x^*, y^*)$ linearly in expectation.

\end{corollary} \label{coro: 1}

% \lizn{Why randomized version? compared with ESGDA?}

\subsection{Convergence under NC-P{\L} condition} \label{sec: 4.3}

%We have proved the convergence of RSGDA for the $\mu$-SCSC function. 
Next, we switch to analyze the property of RSGDA under the NC-P{\L} condition, which is an extension of the result of the $\mu$-SCSC case.  
In general, SGDA and SGDmax can both converge to a Stackelberg-type stationary point in the NC-P{\L} condition. 
However, there is no theoretical analysis ensuring that ESGDA converges in either the NCSC (nonconvex strongly concave) setting or the NC-P{\L} condition~\cite{sebbouh2021randomized}, due to the complicated structure of ESGDA. 
Hence, we analyze the efficiency of RSGDA in the NC-P{\L} condition. 
Our analysis of RSGDA is similar to the analysis of SGDA, which makes the theory simple and transparent.
Specifically, if function $F$ admits the NC-P{\L} condition, we have the following convergence result for Algorithm~\ref{algo: rsgda}.

\begin{theorem} \label{thm: convergence of rsgda}
% Let Assumption \ref{assum: smoothness}, \ref{assum: bounded variance} $\&$ \ref{assum: ncpl condtion} hold. 
Assume that $F$ is $L_1$-smooth and satisfies the NC-P{\L} condition. 
Let $\alpha_k \leq \frac{1}{2L_2}$ and $18\kappa^2 \frac{p}{1-p}  \alpha_k \leq \eta_k \leq \frac{1}{L_1}$ for $k \geq 0$, 
and assume that the stepsizes $\alpha_k$ and $\eta_k$ are square summable but not summable, 
i.e.,
$ \sum_k \alpha_k = + \infty, \sum_k \eta_k = + \infty$, and 
$\sum_k \alpha_k^2 < + \infty, \sum_k \eta_k^2 < + \infty$.
Then for non-increasing $\{ \alpha_k \}$, we have 
\begin{equation}
    \min_{t=0, 1, \dots, k-1} h_t = o \left( \frac{1}{\sum_{j=0}^{k-1} \alpha_j} \right) \to 0, \quad \textrm{almost surely},
\end{equation}
where $h_t = \frac{1}{4} \| \nabla \phi (x_t) \|^2 + \frac{1}{20} \left( \frac{L_1}{\mu} \right)^2 \| \nabla_y F(x_t, y_t) \|^2 + \frac{11}{40} \| \nabla_x F(x_t, y_t) \|^2$.

\end{theorem}

\noindent{\em Remarks.}
% Before we provide a sketch proof of Theorem \ref{thm: convergence of rsgda}, we first discuss what $h_t$ is.
% We claim that 
Here $h_t$ can be viewed as an efficiency measure of the RSGDA.
In fact, we have $\| \nabla F(x_t, y_t) \|^2 \lesssim h_t$ and $\| \nabla \phi (x_t) \|^2 \lesssim h_t$. Thus, we obtain that 
\begin{equation*}
    \min_{t=0, \dots k} \| \nabla F(x_t, y_t) \|^2 + \| \nabla \phi (x_t) \|^2 \lesssim \min_{t=0, \dots, k} h_t,
\end{equation*}
which means that $h_t$ is an upper bound of $\| \nabla F(x_t, y_t) \|^2 + \| \nabla \phi (x_t) \|^2$.
Hence, to obtain a Nash-type stationary point and a Stackelberg-type stationary point simultaneously, it is sufficient to ensure that $\min_{t=0, \dots, k} h_t \to 0$ as $k \to \infty$.
%Hence, roughly speaking, 
In other words, Theorem \ref{thm: convergence of rsgda} essentially states that RSGDA converges to a Nash-type stationary point and a Stackelberg-type stationary point simultaneously under some suitable conditions. 

A detailed proof is provided in Appendix~\ref{app: rsgda}. 
In a nutshell
, to prove Theorem \ref{thm: convergence of rsgda}, inspired by the work \cite{yang2021faster}, we introduce a Lyapunov function 
\begin{equation*} 
V (x, y):= \phi (x) + C (\phi (x) - F (x, y)), 
\end{equation*}
where $C > 0$ is a constant to be determined later. 
Note that for any $(x, y)$, we have $\phi (x) \geq F(x, y)$ according to the definition of $\phi$.
Hence, $V$ is bounded from below.
Now, for any $k \geq 0$, we define
\begin{equation*}
    V_k := V(x_k, y_k) = \phi (x_k) + C(\phi (x_k) - F(x_k, y_k)).
\end{equation*}
Due to the random term $p$ and the stochastic term $z_k$ in the definition of $(x_{k+1}, y_{k+1})$, comparing $V_k$ and $V_{k+1}$ directly is meaningless.
However, from the stochastic process perspective, one can discuss the gap between $\e_k[V_{k+1}]$ and $V_k$, where $\e_k[\cdot]$ is the conditional expectation.
Next, we show that $\{ V_k \}_k$ is similar to a submartingale which means that it is ``non-increasing'' in the conditional expectation meaning.
In particular, we provide an inequality connecting $\e_k[V_{k+1}]$ and $V_k$.
Finally, by applying the Robbins-Siegmund theorem~\cite{robbins1971convergence} to this inequality, we can obtain the almost surely convergence of RSGDA.

% We claim that we can obtain a Stackelberg-type stationary point and a Nash-type stationary point simultaneously from the sequence $\{ x_k, y_k \}_k$ by Theorem \ref{thm: convergence of rsgda}. In fact, on one side, we have $\| \nabla F(x_t, y_t) \|^2 \lesssim h_t$ and $\| \nabla \phi (x_t) \|^2 \lesssim h_t$. Thus, we obtain that 
% \begin{equation*}
%     \min_{t=0, \dots k} \| \nabla F(x_t, y_t) \|^2 + \| \nabla \phi (x_t) \|^2 \lesssim \min_{t=0, \dots, k} h_t.
% \end{equation*}
% On the other side, note that for any given $k$, $\min_{t=0, \dots, k} h_t$ is bounded by $1 / \sum_{j=0}^k \alpha_j$ and $\sum_j \alpha_j = +\infty$. Hence, $\min_{t \geq 0} h_t = 0$ almost surly. Combining the two observations, we get our claim.

Furthermore, we can get the convergence rate of Algorithm \ref{algo: rsgda} from a straightforward observation of Theorem \ref{thm: convergence of rsgda}, which is concluded in the following theorem.

\begin{theorem}[Convergence rate] \label{thm: convergence rate}

Consider the setting of Theorem \ref{thm: convergence of rsgda}, for any $\epsilon > 0$, there are sequences $\{ \alpha_k \}$ and $\{ \eta_k \}$, such that
\begin{equation}
    \min_{t = 0, 1, \dots, k-1} h_t = o \left( k^{- \frac{1}{2} + \epsilon}  \right),
\end{equation}
almost surely.

\end{theorem}

Finally, we analyze two specific versions of RSGDA, i.e., randomized gradient descent ascent (RGDA) and constant step RSGDA.

\noindent\textbf{RGDA}.
We first consider RGDA, where we use the exact gradients in Algorithm \ref{algo: rsgda}.

\begin{corollary} \label{coro: 2}

Let Assumption \ref{assum: smoothness}, \ref{assum: bounded variance} \& \ref{assum: ncpl condtion} hold with $\sigma^2 = 0$. Assume that  $\alpha_k \equiv \alpha$ and $\eta_k \equiv \eta$ for all $k \geq 0$. Moreover, we assume that $\alpha \leq 1 / (2L_2)$ and $18 \frac{p}{1-p} (L_1/\mu)^2 \alpha = \eta \leq 1 / L_1$, where $L_2 = L_1 + \frac{L_1 \kappa}{2}$ and $\kappa = L_1 / \mu$. Then
\begin{equation}
    \min_{t=0,\dots, k} \e [h_t] = O \left( \frac{1}{k} \right).
\end{equation}
    
\end{corollary}

\noindent\textbf{Constant step RSGDA}.
Another variant of RSGDA is choosing constant step sizes for RSGDA. Though RSGDA with constant steps does not converge, we can provide the computation complexity of obtaining an approximate local solution to problem (\ref{problem: stochastic minimax problem}).

\begin{corollary} \label{coro: 3}
    
Let Assumption \ref{assum: smoothness}, \ref{assum: bounded variance} \& \ref{assum: ncpl condtion} hold. Assume that for any $k \geq 0$, we have $\alpha_k \leq 1 / (2L_2)$ and $18 \frac{p}{1-p} (L_1/\mu)^2 \alpha_k = \eta_k \leq 1 / L_1$, where $L_2 = L_1 + \frac{L_1 \kappa}{2}$ and $\kappa = L_1 / \mu$. Moreover, assume that $\alpha_k \equiv \alpha$ and $\eta_k \equiv \eta$ for all $k \geq 0$. Then for any $\epsilon > 0$, if $k = O (\epsilon^2)$, then
\begin{equation}
    \min_{t=0, \dots, k} \e [h_t] \leq \epsilon.
\end{equation}
    
\end{corollary}

%\begin{remark}
\noindent{\em Remarks}.
Corollary \ref{coro: 3} indicates that the computation complexity of RSGDA is the same as SGDA. 
In other words, RSGDA is as fast as SGDA.
In fact, for any $\epsilon > 0$, SGDA (see like \citet{yang2021faster}) provides a point $x$ such that $\| \nabla \phi (x) \|^2 \leq \epsilon$ in $O(1/\epsilon^2)$ steps.
On the other side, note that $h_k \geq\frac{1}{4} \| \phi (x_k) \|^2$ for any $k \geq 0$.
Hence, Corollary \ref{coro: 3} states that RSGDA provides a point $x$ such that $\| \nabla \phi (x) \|^2 \leq \epsilon$ in $O(1/\epsilon^2)$ steps, which coincides with the results of SGDA.
It means that RSGDA shares the same theoretical computation complexity as SGDA.
%\yy{this paragraph seems to be a mismatch to the above clrollary.}

%\end{remark}

%\yy{overall, the above two sections are still messy to me. consider to add a figure showing the logic flow (this might be rare but I guess a figure may help the readers to follow), and a table to summarize the assumptions/results/etc. of existing methods and your method.}

\subsection{Selection of p} \label{sec: 4.4}

In the last part, we discuss about the value $p$ in applications.
To the end of this subsection, we always assume that the stepsizes are constants.
Theoretically, for any $p \in (0, 1)$, RSGDA converges under some suitable choices of stepsizes. 
However, the empirical experiments show that RSGDA performs better in some $p$ than others.
We propose an intuitive idea to determine the value $p$ in this part.

Recall that $h_k$ is introduced to measure the efficiency of RSGDA.
Thus, for any fixed $n$, the smaller value $\sum_{k=1}^n h_k$, the faster RSGDA converges numerically.
Hence, the basic idea of determining $p$ is to minimize the upper bound of $\sum_{k=1}^n h_k$ by choosing suitable parameter $p$.
Specifically, we choose $p$ according to the following rule.

\begin{proposition} \label{prop: 1}

Consider the setting of Corollary \ref{coro: 3}.
For any initial point $(x_0, y_0)$, There are constant $M_1, M_2 > 0$ independent of $\{ z_1, x_1, y_1, \dots \}$ (the explicit form of $M$ is given in the Appendix), such that the optimal $p$ is given as follow.
\begin{equation}
p =
\begin{dcases}
    \min \left\{ \frac{M_1}{\sqrt{n}} - \frac{M_2}{n} ,\frac{L_2}{9 L_1 \kappa^2 + L_2} \right\}, & \quad \sigma > 0; \\
    \frac{L_2}{9 L_1 \kappa^2 + L_2} & \quad \sigma = 0.
\end{dcases}
\end{equation}
where $n$ is the step number.

\end{proposition}

We give a proof sketch here.
Throughout a complicated computation, one can show that 
\begin{equation} \label{equ: 4.11}
    \sum_{k=1}^n h_k \leq \frac{1}{\alpha p} (V_1 - \inf V) + \frac{18^2 p n}{2 (1-p)} \kappa^4 L_1 \alpha \sigma^2 + \mathrm{const},
\end{equation}
where $V_1$ and $\inf V$ are two constants determined by the initial values, and $\mathrm{const}$ is a constant.
We aim to find the minimization of RHS of \eqref{equ: 4.11} with respect to $p$.
The basic theory of calculus shows that the minimization of RHS of \eqref{equ: 4.11} can be obtained at the value offered in Proposition \ref{prop: 1}.

Roughly speaking, Proposition \ref{prop: 1} states that one should choose the probability $p$ obeying the following rules.
In case the estimation of the gradient is exact, i.e., there is no variance, the hyperparameter $p$ should be as large as possible.
In case $\sigma > 0$, $p$ is the smaller one between $O(1 / \sqrt{n})$ and a constant $\frac{L_2}{9 L_1 \kappa^2 + L_2}$.
Hence, if $n$ is small, $p$ should be a constant; if $n$ is sufficiently large, $p$ should decrease in the rate $O(1/ \sqrt{n})$.
In other words, RSGDA prefers a large probability $p$ when the iteration number $n$ is small and a small $p$ with order $O(1/\sqrt{n})$ in the case $n$ is large.
Hence, we propose an intuitive method called Adaptive-RSGDA (AdaRSGDA) to adjust the parameter $p$.
First, we choose an integer $N_1 > 0$, a constant $N_2 > 0$, and initial probability $p_0$.
For the iteration step $n < N_1$, we fix $p = p_0$.
For the case $n \geq N_1$, we choose $p = 1/([(n-N_1) / N_2] +1)$ for each $N_2$ steps, where $[x]$ denotes the largest integer smaller than $x$.
We use this simple notation to estimate the term $O(1/\sqrt{n})$.
A straightforward observation shows that AdaRSGDA performs like SGDA at the starting $N_1$ steps and naturally changes to ESGDA from the $N_1 + 1$ step.

% In summary, an intuitive choice of $p$ is that for the first several iterations, $p$ should be large enough, then $p$ decreases gradually.
% In other words, RSGDA should take more gradient descent steps at the beginning of the iteration and gradually takes more gradient ascent steps.

\section{Numerical experiments} \label{sec: experiment}
We conduct experiments on both synthetic data and real data.\footnote{{The code, together with the corresponding dataset, are uploaded to \url{https://figshare.com/s/207ab3663e2acd28d8dd}.}} The experiments are mainly designed to answer the following questions:
\begin{itemize}
    \item Does RSGDA perform consistently with ESGDA?
    %\item Is RSGDA as fast as SGDA as indicated by the theoretical result?
    \item Can RSGDA outperform ESGDA under certain circumstances?
    \item How efficient is the proposed selection strategy for probability $p$ in RSGDA?
\end{itemize}

\subsection{Experiments on synthetic data}

\subsubsection{MLP Wasserstein GANs.} 
Following the problem setting of~\citet{loizou2020stochastic}, we first use the WGAN to approximate a given multivariate Gaussian distribution.
We consider the WGAN task in which the discriminator and generator are both modeled with MLPs~\citep{lei2020sgd}. 
%Following the problem setting of~\citet{loizou2020stochastic}, 
%we use WGAN to approximate a given multivariate Gaussian distribution. 
Elaborately, the real data is drawn from a normal distribution $\mathcal{N}(\mu^*, \sigma^*)$ where $\mu^*=(0.5; -1.5)$ and $\sigma^*=(0.1; 0.3)$, and the fake data is generated by a generator denoted by $g_{\theta}(z)$, where $z$ is drawn from the standard Gaussian distribution.
The discriminator is defined as $f_{w}(x)$, where input $x$ can be either a real example or a fake example. 
The minimax problem of this task can be formulated as
\begin{equation*}
\min_{\theta} \max_{w} F(w, \theta) := ~\e_{x, z} [f_{w}(x) - f_{w}(g_{\theta}(z))]. 
\end{equation*}
In this experiment, we fix the batch size to $100$, learning rates $\alpha=0.01$ and $\eta=0.01$ for all approaches. 
Each reported result in the following is the average of $5$ repeated experiments.

\begin{figure}[tbp]
\centering  
\subfigure[$k = 1$]{
\label{fig: 1a}
\includegraphics[width=0.455\columnwidth]{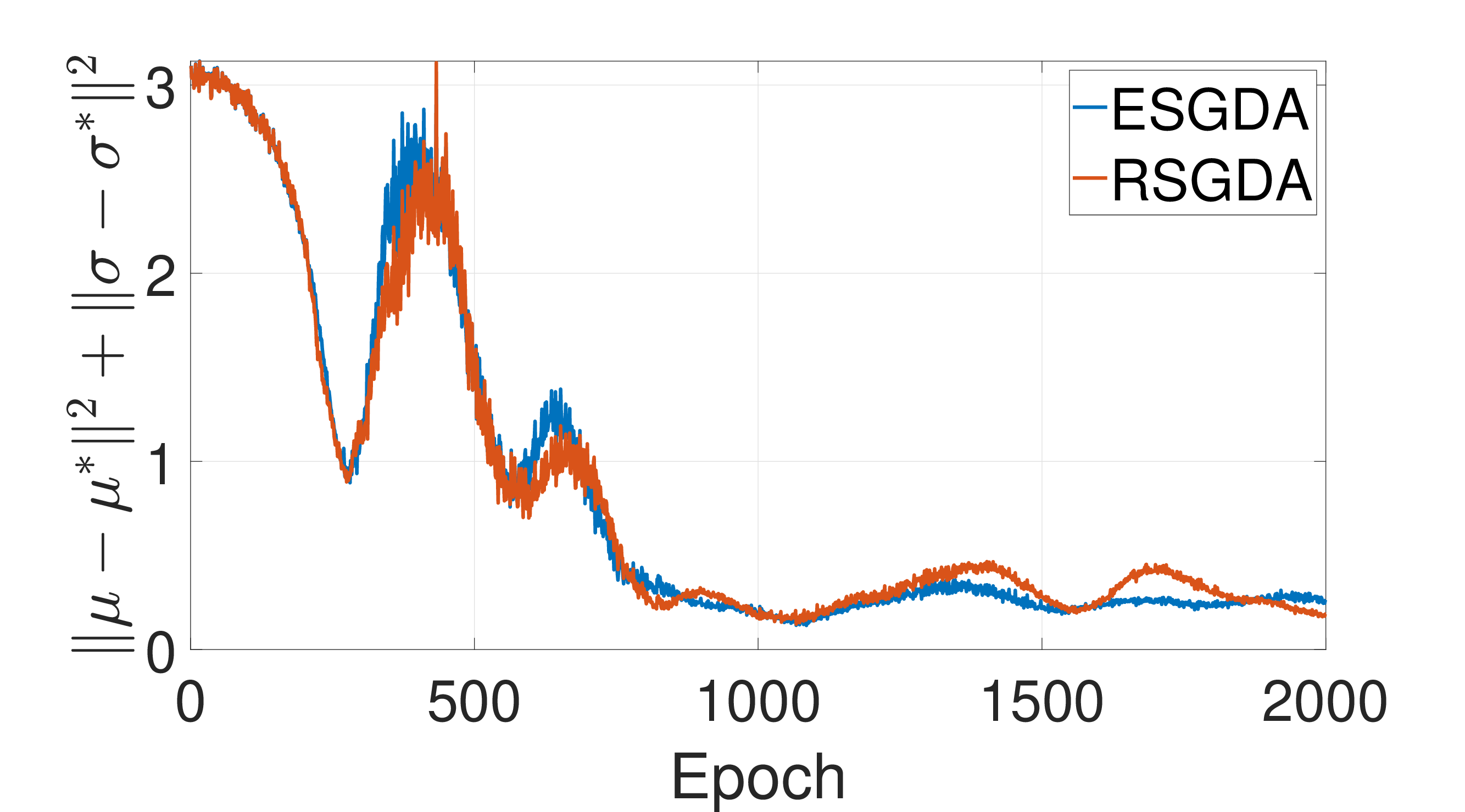}}
\subfigure[$k=3$]{
\label{fig: 1b}
\includegraphics[width=0.455\columnwidth]{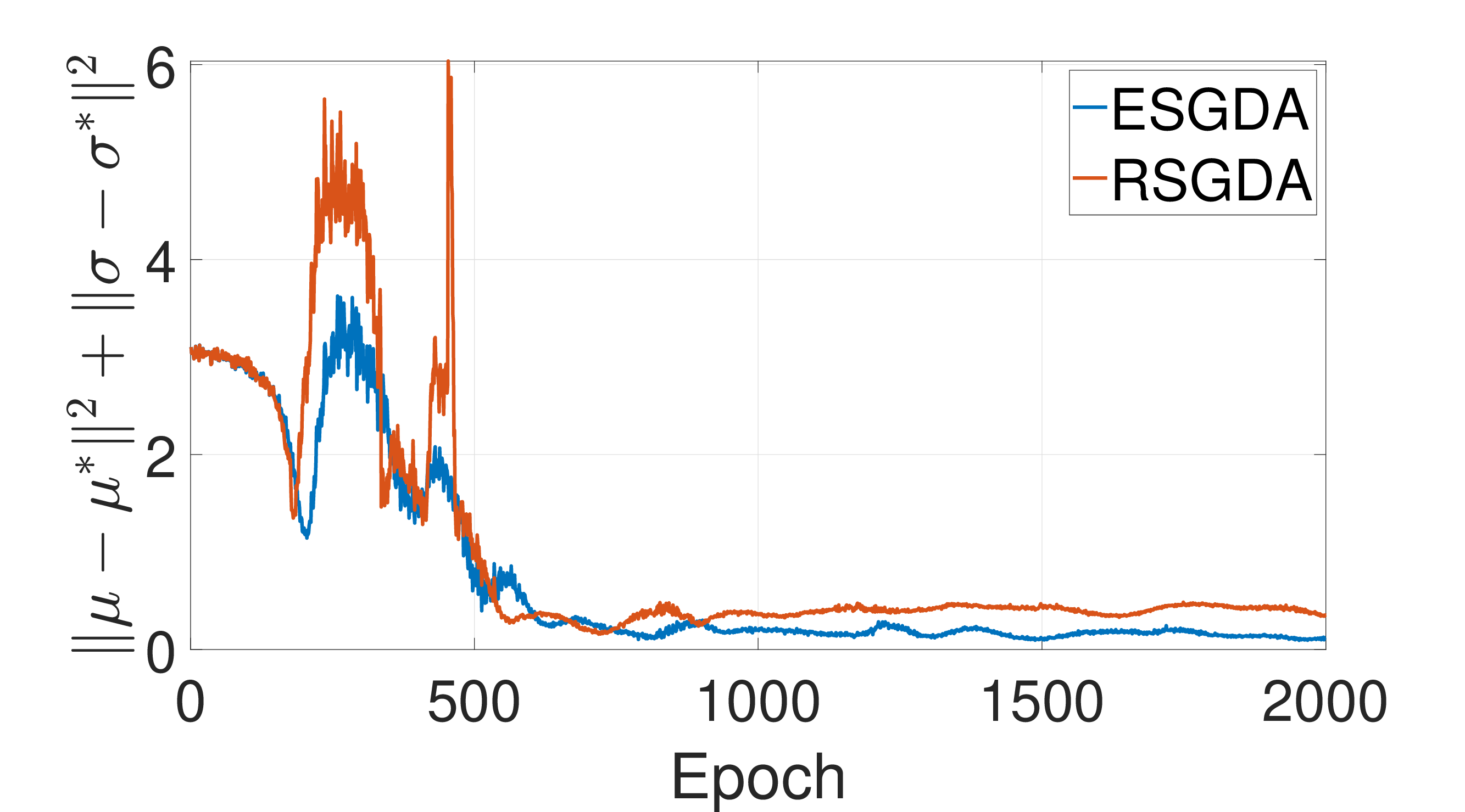}} \\
\subfigure[$k = 5$]{
\label{fig: 1c}
\includegraphics[width=0.455\columnwidth]{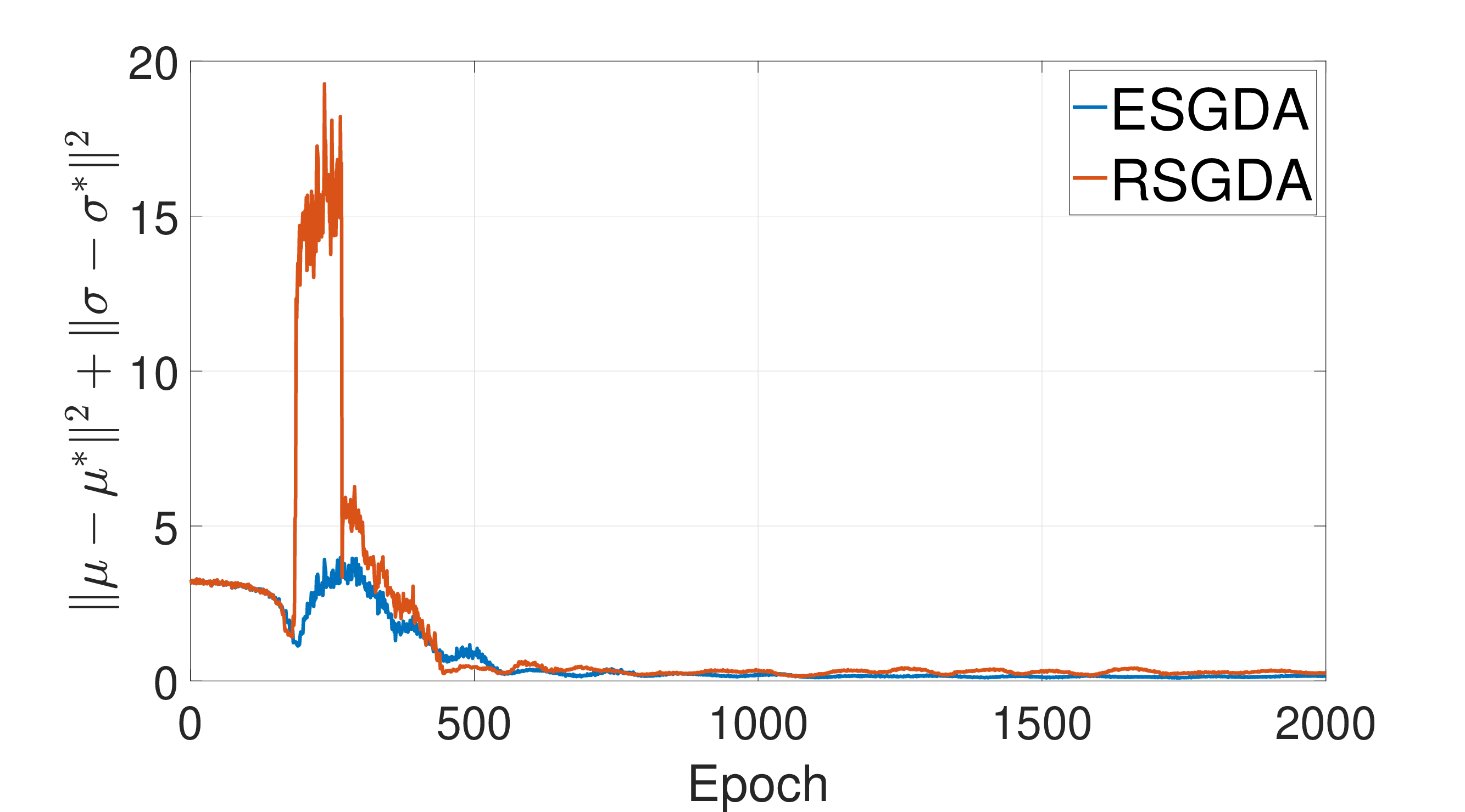}}
\subfigure[$k = 7$]{
\label{fig: 1d}
\includegraphics[width=0.455\columnwidth]{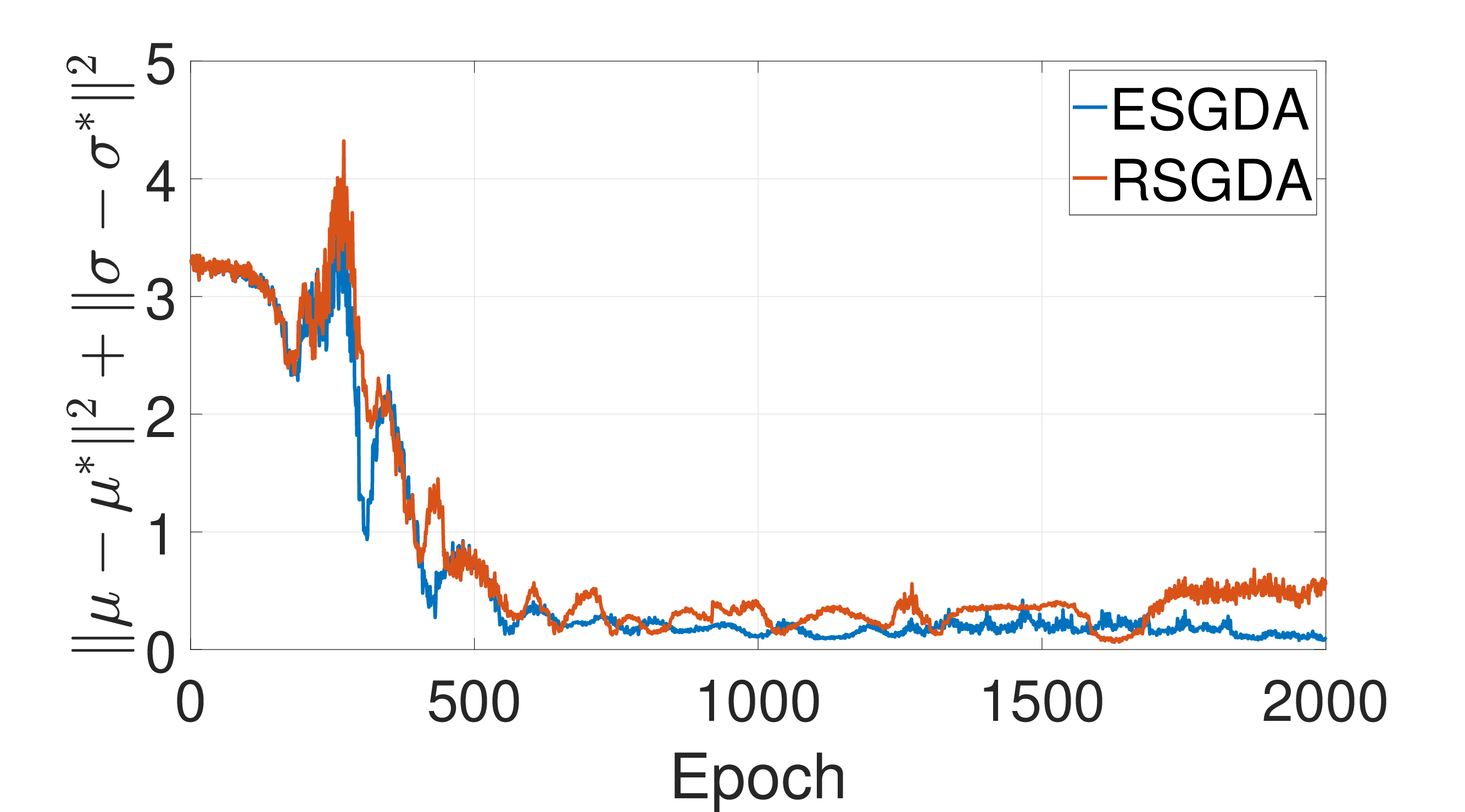}}
\caption{The distances to the optimal solution of ESGDA and RSGDA with different $(m,p)$ settings on WGAN. RSGDA and ESGDA generally have consistent performance.
%\yy{the words in the figures are too small. btw, what it the y-axis. it is the same with figure 2, right? then why two difference names.}
}
\label{fig: 0}
\end{figure}

To confirm the consistency between ESGDA and RSGDA, we choose different sizes of the inner loop in ESGDA, and set the probability $p$ in RSGDA accordingly.
Specifically, we set pairs $(m,p)$ by $\{(1, \frac{1}{2}), (3, \frac{1}{4}), (5, \frac{1}{6}), (7, \frac{1}{8}) \}$. 
The training curves of ESGDA and RSGDA are shown in Figure~\ref{fig: 0}, where the y-axis measures the distance to the optimal solution $(\mu^*, \sigma^*)$.
We can observe that, except for a few outliers caused by randomness, the behavior of RSGDA is consistent to that of ESGDA. 

\begin{figure}[tb]
    \centering
    \includegraphics[width=1.0\columnwidth]{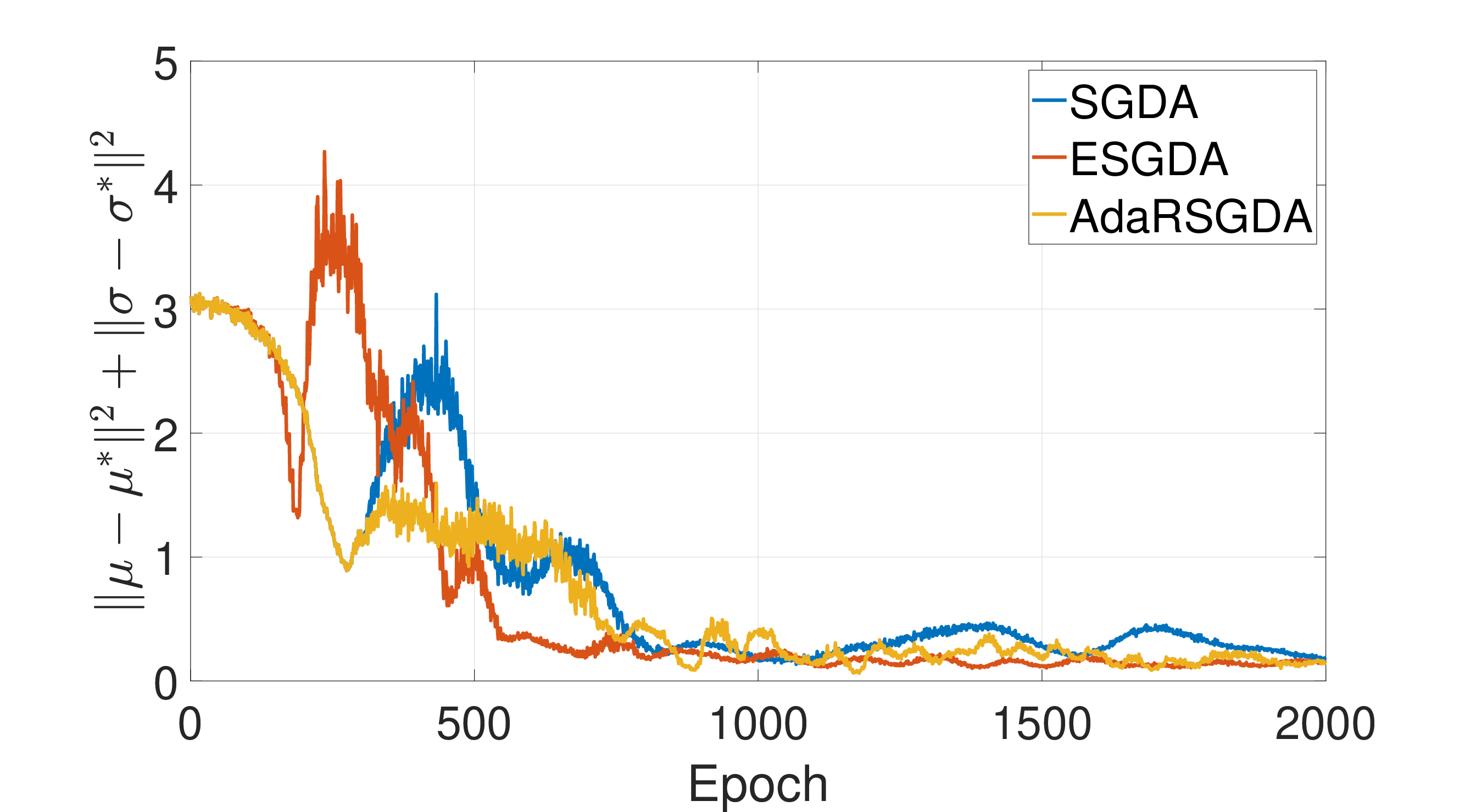}
    \caption{The distances to the optimal solution of  SGDA, ESGDA, and AdaRSGDA on WGAN. AdaRSGDA performs closely to the counterparts while it is more stable.}
    \label{fig: 1}
\end{figure}

We also evaluate the efficiency of our method AdaRSGDA.
Figure~\ref{fig: 1} provides a comparison of AdaRSGDA to ESGDA and SGDA. %by measuring the distance to the optimal solution $(\mu^*, \sigma^*)$.
For ESGDA, we set $m=5$, 
and for AdaRSGDA, we simply define $N_1 = N_2 = 300$. 
We can first observe that AdaRSGDA generally performs better than SGDA. 
Also, AdaRSGDA is more stable and converges faster than ESGDA at the beginning of the iterations. 
Moreover, as $n$ grows larger, AdaRSGDA has a refined estimation of the optimal solution which is as good as ESGDA.

\subsubsection{Robust non-linear regression.}

We next consider the robust non-linear regression problem proposed by~\citet{yang2021faster}. %to further evaluate the effect of hyperparameter $p$ on the convergence rate.
We generate a dataset consisting of $1000$ data points in $500$ dimensions, sampled from the normal distribution $\mathcal{N}(0, 1)$. %where 
% \lizn{sampled from the Gaussian distribution with mean $0$ and variance $1$.}
The target value $y$ is sampled by a random linear model with an additional noise. 
We define $f_w (z)$ as an MLP model with the parameter $w$.
The goal of robust non-linear regression model is to solve the following minimax problem:
\begin{equation*}
    \min_{w} \max_{y} ~\frac{1}{n} \sum_{j=1}^n \frac{1}{2} \| f_w (x_i) - y \|^2 - \frac{1}{2} \| y - y_i \|^2.
\end{equation*}

In this experiment, we mainly focus on the effect of hyperparameter $p$ on the convergence rate. 
%We confirm the effect of hyperparameter $p$ on the convergence rate of RSGDA and compare the efficiency of RSGDA with ESGDA and SGDA. 
We set the batch size to $1000$, learning rate $\alpha=5e-4$ and $\eta=5$. 
We set the probability $p$ of RSGDA by $p=0.2$ and $p = 0.8$, and compare the convergence rate of RSGDA with ESGDA and SGDA. The size of the inner loop is set as $m=4$ for ESGDA. 
Note that the choice $p=0.8$ leads to more update steps in variable $x$ than $y$.
Meanwhile, ESGDA always takes more update steps in $y$ than $x$.
% Hence, RSGDA with parameter $p=0.8$ goes beyond the ESGDA framework.

\begin{figure}[tb]
    \centering
    \includegraphics[width=1.0\columnwidth]{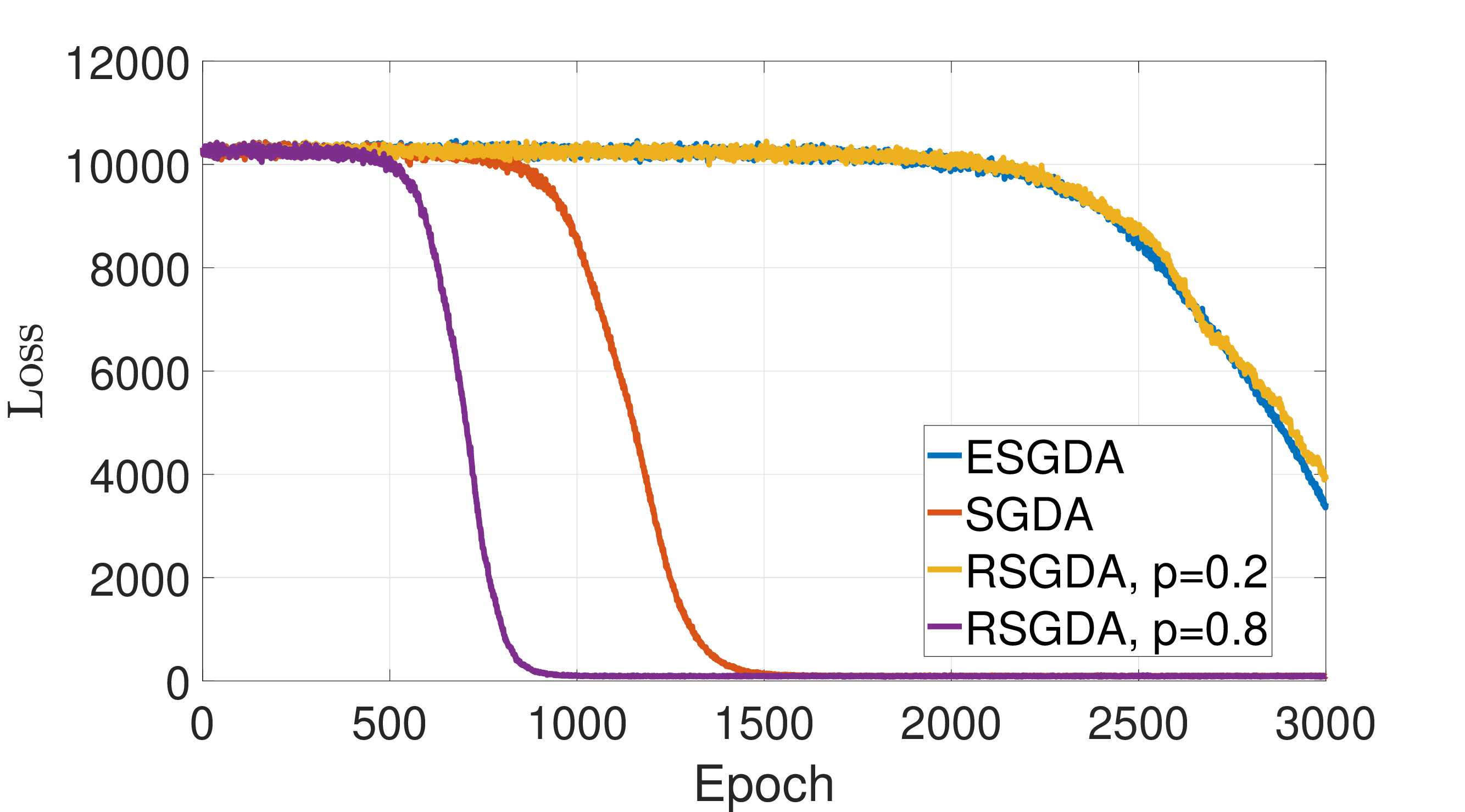}
    \caption{The loss curves of ESGDA, SGDA, and RSGDA on non-linear regression. RSGDA with large $p$ broadens the applicable range and outperforms ESGDA and SGDA.}
    \label{fig: regression}
\end{figure}

The loss curves are shown in Figure~\ref{fig: regression}.
There are some remarkable observations.
First, the curve of RSGDA with $p=0.2$ coincides with that of ESGDA, further indicating the consistency between RSGDA and ESGDA ($p=0.2$ is consistent with $m=4$). % which shows the consistency of ESGDA and RSGDA.
Second, we observe that RSGDA could converge even faster than SGDA, which is consistent with our analysis that RSGDA converges as fast as SGDA.
Third, with larger $p$, RSGDA also converges faster than ESGDA. For example, RSGDA with $p=0.8$ converges fastest in the four curves, which indicates that it is better to take more outer steps in $x$ rather than inner steps in $y$ in this model.
In other words, ESGDA is not always the best choice, and RSGDA can broaden the applicable range by choosing a large parameter $p$ (e.g., $p > 0.5$). % it is meaningful to discuss the performance of RSGDA with parameter $p > 0.5$.

\subsection{Experiments on real data}
\subsubsection{Adversarial training.}
Finally, we study the adversarial training task with real data. Adversarial training aims to ensure the model to be robust against adversarial perturbations. 
Given the training data distribution $\mathcal{D}$, and letting $f_{\bm{w}}(\cdot)$ denote the classifier parameterized by $w$, and $L(\cdot, \cdot)$ denote the cross-entropy loss, the minimax problem of adversarial training can be formulated as~\citep{sinha2017certifying}
\begin{equation*}
    \min_{\bm{w}} \max_{\bm{\delta} \in S_p} ~\e_{(\bm{x}, \bm{y}) \sim \mathcal{D}} L(f_{\bm{w}}(\bm{x}+\bm{\delta}), \bm{y}),
\end{equation*}
where $S_p$ is the $\ell_p$ norm ball introduced to make the perturbation $\bm{\delta}$ small enough in the sense of $\ell_p$ norm. 

\begin{figure}[t]
    \centering
    \includegraphics[width=1.0\columnwidth]{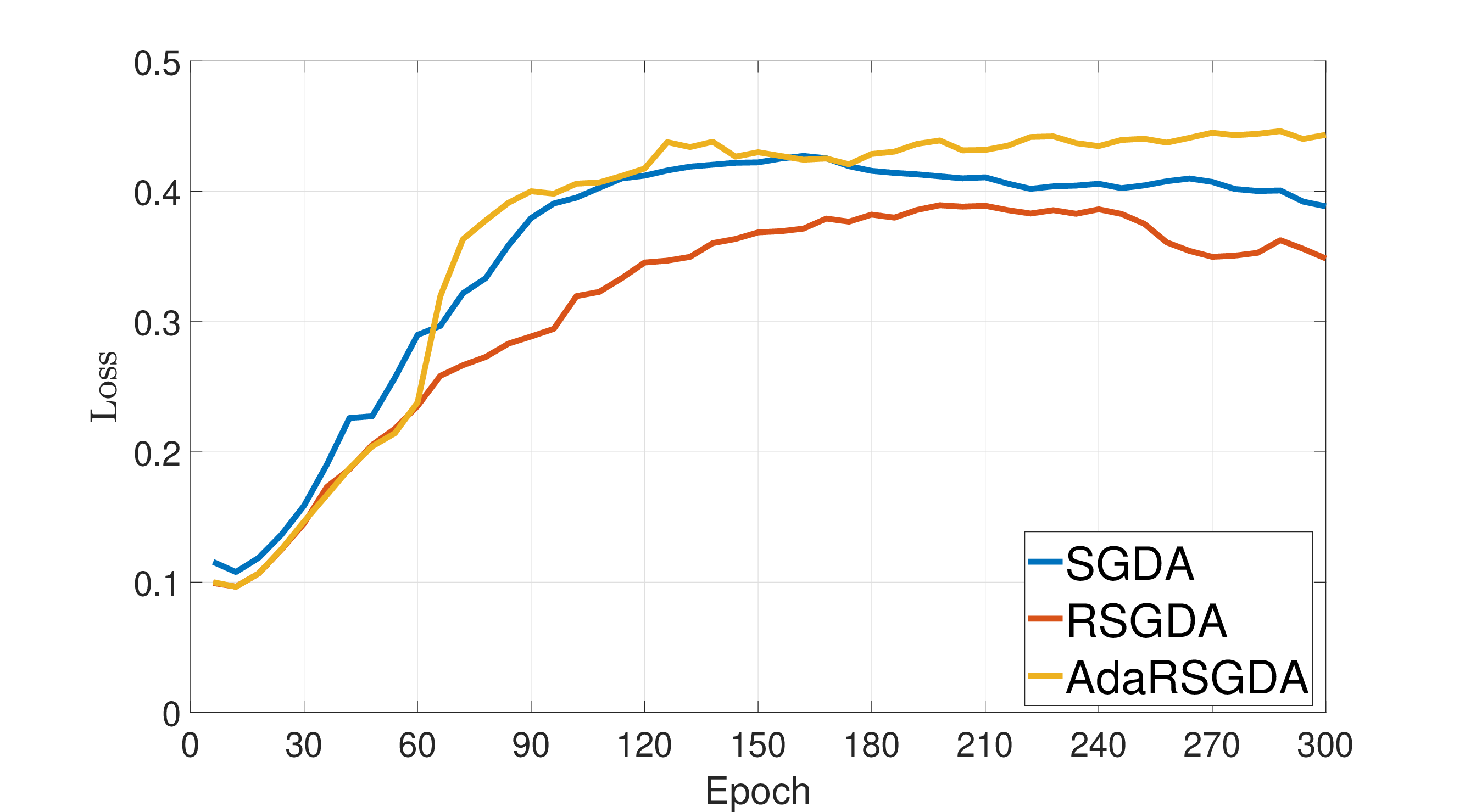}
    \caption{The training loss curves of SGDA, RSGDA, and AdaRSGDA on adversarial training.}
    \label{fig: adv loss}
\end{figure}

We conduct the experiment on the MNIST dataset~\citep{lecun2010mnist}. 
In this task, $f_{\bm{w}}$ is specific to the LeNet-5 model~\citep{lecun1998gradient}, and $S_p$ is defined as $S_p := \{\bm{\delta} \mid \|\bm{\delta}\|_{\infty} \leq 0.3\}$~\citep{madry2017towards}. 
The step sizes $\alpha$ and $\eta$ of the gradient descent and ascent are both fixed to $0.1$. 
%We evaluate the efficiency of AdaRSGDA with the following setting. 
For RSGDA, we set $p = 0.5$.
For AdaRSGDA, we set  $p_0=0.5$ and $N_1=N_2=60$.  %inner loop size $m=5$ for ESGDA,
The training loss curve is plotted in Figure \ref{fig: adv loss} and the accuracy is given in Table \ref{tab: accuracy}. 
% \yy{is it RSGDA or SGDA. the figure and the table are not consistent}

% To confirm the consistency between ESGDA and RSGDA, we choose different sizes of the inner loop in ESGDA, and set the probability $p$ in RSGDA accordingly.
% For details, we set pairs $(k,p)$ by 
% \begin{equation*}
% (k, p) \in \left\{(1, \frac{1}{2}), (3, \frac{1}{4}), (5, \frac{1}{6}), (7, \frac{1}{8}) \right\}.
% \end{equation*}

\begin{table}[tb]
\begin{center}
\caption{The accuracy ($\%$) on MNIST dataset. AdaRSGDA performs better than the plain SGDA and RSGDA on adversarial data.}\label{tab: accuracy}
\begin{threeparttable}
\begin{tabular}{c|c|c}
\toprule 
{\textbf{Methods}} & {\textbf{Benign data}} & {\textbf{Adversarial data}}  \\
\midrule
%ESGDA & $98.2$ & $88.7$ \\
RSGDA & $99.1$ & $57.3$ \\
SGDA & $98.9$ & $75.4$ \\
AdaRSGDA & $98.9$ & $83.7$ \\
\bottomrule
\end{tabular}
\end{threeparttable}
\end{center}
\end{table}

Figure \ref{fig: adv loss} shows that the algorithms perform well on the training data, and the losses are close.
The accuracy in Table \ref{tab: accuracy} also supports this conclusion.
Table \ref{tab: accuracy} shows that RSGDA, SGDA, and AdaRSGDA are very close and all achieve high accuracy on the benign data.
However, their performance deviates on the adversarial data, and RSGDA and SGDA perform much worse than AdaRSGDA. 
%That is, SGDA performs much worse than ESGDA and AdaRSGDA on the adversarial (adv) data. 
RSGDA and SGDA only achieve $57.3\%$ and $75.4\%$ accuracies on the adversarial data, while AdaRSGDA has $83.7\%$ accuracy. %ESGDA has $88.7\%$ accuracy and 
Hence, training the same objective function with AdaRSGDA is more robust to adversarial attacks than SGDA and RSGDA in this real data application. The reasons behind this difference need further investigation and we leave them as future work.
A possible reason is that the lack of delay of $p$ slows down the convergence of SGDA and RSGDA.
%\yy{do we have reasons for this?}

%\yy{somewhere we may need to re-mention things  we have said in the intro, eg. RSGDA is as fast as SGDA.}

\section{Conclusion}

In this paper, we establish a new convergence analysis for RSGDA under milder assumptions beyond strong concavity.
In particular, we show that RSGDA has the same convergence rate with SGDA in the NC-P{\L} condition, which is a better complexity compared with the previous analysis of RSGDA.
Moreover, we propose an intuitive method to choose the parameter $p$ for RSGDA and confirm its efficiency in experiments. 
There are quite a few works worthy to be considered in the future. 
On the one hand, it is still an open question whether there exists a tighter complexity bound for the minimax problems in the NC-P{\L} condition or even more general cases. 
On the other hand, a more efficient adaptive RSGDA for minimax optimization needs further investigation. 
% First, how to design adaptive algorithms for the minimax problems?
% Moreover, it is still an open question that what is the lower complexity bound for the minimax problems in the NC-P{\L} condition.
% Furthermore, does RSGDA still hold for the general minimax problems?

\bibliography{src/reference}

\begin{thebibliography}{43}
\providecommand{\natexlab}[1]{#1}

\bibitem[{Bolte, Daniilidis, and Lewis(2006)}]{bolte2007lojasiewicz}
Bolte, J.; Daniilidis, A.; and Lewis, A. 2006.
\newblock The {{\L}}ojasiewicz inequality for nonsmooth subanalytic functions
  with applications to subgradient dynamical systems.
\newblock \emph{SIAM J. Optim.}, 17(4): 1205--1223.

\bibitem[{Bolte et~al.(2007)Bolte, Daniilidis, Lewis, and
  Shiota}]{bolte2007clarke}
Bolte, J.; Daniilidis, A.; Lewis, A.; and Shiota, M. 2007.
\newblock Clarke subgradients of stratifiable functions.
\newblock \emph{SIAM J. Optim.}, 18(2): 556--572.

\bibitem[{Bo{\c{t}} and B{\"o}hm(2020)}]{boct2020alternating}
Bo{\c{t}}, R.~I.; and B{\"o}hm, A. 2020.
\newblock Alternating proximal-gradient steps for (stochastic)
  nonconvex-concave minimax problems.
\newblock \emph{arXiv preprint arXiv:2007.13605}.

\bibitem[{Cesa-Bianchi and Lugosi(2006)}]{cesa2006prediction}
Cesa-Bianchi, N.; and Lugosi, G. 2006.
\newblock \emph{Prediction, learning, and games}.
\newblock Cambridge university press.

\bibitem[{Chen et~al.(2021)Chen, Zhou, Xu, and Liang}]{chen2021proximal}
Chen, Z.; Zhou, Y.; Xu, T.; and Liang, Y. 2021.
\newblock Proximal Gradient Descent-Ascent: Variable Convergence under K {\L}
  Geometry.
\newblock \emph{arXiv preprint arXiv:2102.04653}.

\bibitem[{Dai et~al.(2018)Dai, Shaw, Li, Xiao, He, Liu, Chen, and
  Song}]{dai2018sbeed}
Dai, B.; Shaw, A.; Li, L.; Xiao, L.; He, N.; Liu, Z.; Chen, J.; and Song, L.
  2018.
\newblock Sbeed: Convergent reinforcement learning with nonlinear function
  approximation.
\newblock In \emph{International Conference on Machine Learning}, 1125--1134.
  PMLR.

\bibitem[{Daskalakis and Panageas(2018)}]{daskalakis2018limit}
Daskalakis, C.; and Panageas, I. 2018.
\newblock The limit points of (optimistic) gradient descent in min-max
  optimization.
\newblock \emph{Advances in neural information processing systems}, 31.

\bibitem[{Diakonikolas, Daskalakis, and
  Jordan(2021)}]{diakonikolas_2021_efficient}
Diakonikolas, J.; Daskalakis, C.; and Jordan, M. 2021.
\newblock Efficient Methods for Structured Nonconvex-Nonconcave Min-Max
  Optimization.
\newblock In Banerjee, A.; and Fukumizu, K., eds., \emph{Proceedings of The
  24th International Conference on Artificial Intelligence and Statistics},
  volume 130 of \emph{Proceedings of Machine Learning Research}, 2746--2754.
  PMLR.

\bibitem[{Farnia and Ozdaglar(2021)}]{farnia2021train}
Farnia, F.; and Ozdaglar, A. 2021.
\newblock Train simultaneously, generalize better: Stability of gradient-based
  minimax learners.
\newblock In \emph{International Conference on Machine Learning}, 3174--3185.
  PMLR.

\bibitem[{Fiez et~al.(2021)Fiez, Ratliff, Mazumdar, Faulkner, and
  Narang}]{fiez2021global}
Fiez, T.; Ratliff, L.; Mazumdar, E.; Faulkner, E.; and Narang, A. 2021.
\newblock Global convergence to local minmax equilibrium in classes of
  nonconvex zero-sum games.
\newblock \emph{Advances in Neural Information Processing Systems}, 34:
  29049--29063.

\bibitem[{Goodfellow et~al.(2014)Goodfellow, Pouget-Abadie, Mirza, Xu,
  Warde-Farley, Ozair, Courville, and Bengio}]{goodfellow2014generative}
Goodfellow, I.; Pouget-Abadie, J.; Mirza, M.; Xu, B.; Warde-Farley, D.; Ozair,
  S.; Courville, A.; and Bengio, Y. 2014.
\newblock Generative adversarial nets.
\newblock \emph{Advances in neural information processing systems}, 27.

\bibitem[{Goodfellow, Shlens, and Szegedy(2014)}]{goodfellow2014explaining}
Goodfellow, I.~J.; Shlens, J.; and Szegedy, C. 2014.
\newblock Explaining and harnessing adversarial examples.
\newblock \emph{arXiv preprint arXiv:1412.6572}.

\bibitem[{Heusel et~al.(2017)Heusel, Ramsauer, Unterthiner, Nessler, and
  Hochreiter}]{heusel2017gans}
Heusel, M.; Ramsauer, H.; Unterthiner, T.; Nessler, B.; and Hochreiter, S.
  2017.
\newblock Gans trained by a two time-scale update rule converge to a local nash
  equilibrium.
\newblock \emph{Advances in neural information processing systems}, 30.

\bibitem[{Jin, Netrapalli, and Jordan(2020)}]{jin2020local}
Jin, C.; Netrapalli, P.; and Jordan, M. 2020.
\newblock What is local optimality in nonconvex-nonconcave minimax
  optimization?
\newblock In \emph{International Conference on Machine Learning}, 4880--4889.
  PMLR.

\bibitem[{Karimi, Nutini, and Schmidt(2016)}]{karimi2016linear}
Karimi, H.; Nutini, J.; and Schmidt, M. 2016.
\newblock Linear convergence of gradient and proximal-gradient methods under
  the polyak-{{\L}}ojasiewicz condition.
\newblock In \emph{Joint European Conference on Machine Learning and Knowledge
  Discovery in Databases}, 795--811. Springer.

\bibitem[{Kurdyka(1998)}]{kurdyka1998gradients}
Kurdyka, K. 1998.
\newblock On gradients of functions definable in o-minimal structures.
\newblock \emph{Ann. Inst. Fourier (Grenoble)}, 48(3): 769--783.

\bibitem[{LeCun et~al.(1998)LeCun, Bottou, Bengio, and
  Haffner}]{lecun1998gradient}
LeCun, Y.; Bottou, L.; Bengio, Y.; and Haffner, P. 1998.
\newblock Gradient-based learning applied to document recognition.
\newblock \emph{Proceedings of the IEEE}, 86(11): 2278--2324.

\bibitem[{LeCun, Cortes, and Burges(2010)}]{lecun2010mnist}
LeCun, Y.; Cortes, C.; and Burges, C. 2010.
\newblock MNIST handwritten digit database.
\newblock \emph{ATT Labs [Online]. Available:
  http://yann.lecun.com/exdb/mnist}, 2.

\bibitem[{Lee and Kim(2021)}]{lee2021fast}
Lee, S.; and Kim, D. 2021.
\newblock Fast Extra Gradient Methods for Smooth Structured
  Nonconvex-Nonconcave Minimax Problems.
\newblock In Beygelzimer, A.; Dauphin, Y.; Liang, P.; and Vaughan, J.~W., eds.,
  \emph{Advances in Neural Information Processing Systems}.

\bibitem[{Lei et~al.(2020)Lei, Lee, Dimakis, and Daskalakis}]{lei2020sgd}
Lei, Q.; Lee, J.; Dimakis, A.; and Daskalakis, C. 2020.
\newblock Sgd learns one-layer networks in wgans.
\newblock In \emph{International Conference on Machine Learning}, 5799--5808.
  PMLR.

\bibitem[{Li et~al.(2022)Li, Farnia, Das, and Jadbabaie}]{li_2022_convergence}
Li, H.; Farnia, F.; Das, S.; and Jadbabaie, A. 2022.
\newblock On Convergence of Gradient Descent Ascent: A Tight Local Analysis.
\newblock In Chaudhuri, K.; Jegelka, S.; Song, L.; Szepesvari, C.; Niu, G.; and
  Sabato, S., eds., \emph{Proceedings of the 39th International Conference on
  Machine Learning}, volume 162 of \emph{Proceedings of Machine Learning
  Research}, 12717--12740. PMLR.

\bibitem[{Lin, Jin, and Jordan(2020{\natexlab{a}})}]{lin2020gradient}
Lin, T.; Jin, C.; and Jordan, M. 2020{\natexlab{a}}.
\newblock On gradient descent ascent for nonconvex-concave minimax problems.
\newblock In \emph{International Conference on Machine Learning}, 6083--6093.
  PMLR.

\bibitem[{Lin, Jin, and Jordan(2020{\natexlab{b}})}]{lin2020near}
Lin, T.; Jin, C.; and Jordan, M.~I. 2020{\natexlab{b}}.
\newblock Near-optimal algorithms for minimax optimization.
\newblock In \emph{Conference on Learning Theory}, 2738--2779. PMLR.

\bibitem[{Loizou et~al.(2020)Loizou, Berard, Jolicoeur-Martineau, Vincent,
  Lacoste-Julien, and Mitliagkas}]{loizou2020stochastic}
Loizou, N.; Berard, H.; Jolicoeur-Martineau, A.; Vincent, P.; Lacoste-Julien,
  S.; and Mitliagkas, I. 2020.
\newblock Stochastic hamiltonian gradient methods for smooth games.
\newblock In \emph{International Conference on Machine Learning}, 6370--6381.
  PMLR.

\bibitem[{{\L}ojasiewicz(1963)}]{lojasiewicz1963propriete}
{\L}ojasiewicz, S. 1963.
\newblock Une propri{\'e}t{\'e} topologique des sous-ensembles analytiques
  r{\'e}els.
\newblock In \emph{Les {\'E}quations aux {D}{\'e}riv{\'e}es {P}artielles
  ({P}aris, 1962)}, 87--89. {\'E}ditions du Centre National de la Recherche
  Scientifique (CNRS).

\bibitem[{Luo et~al.(2020)Luo, Ye, Huang, and Zhang}]{luo_2020_advances}
Luo, L.; Ye, H.; Huang, Z.; and Zhang, T. 2020.
\newblock Stochastic Recursive Gradient Descent Ascent for Stochastic
  Nonconvex-Strongly-Concave Minimax Problems.
\newblock In Larochelle, H.; Ranzato, M.; Hadsell, R.; Balcan, M.; and Lin, H.,
  eds., \emph{Advances in Neural Information Processing Systems}, volume~33,
  20566--20577. Curran Associates, Inc.

\bibitem[{Madry et~al.(2017)Madry, Makelov, Schmidt, Tsipras, and
  Vladu}]{madry2017towards}
Madry, A.; Makelov, A.; Schmidt, L.; Tsipras, D.; and Vladu, A. 2017.
\newblock Towards deep learning models resistant to adversarial attacks.
\newblock \emph{arXiv preprint arXiv:1706.06083}.

\bibitem[{Mazumdar, Ratliff, and Sastry(2020)}]{mazumdar2020gradient}
Mazumdar, E.; Ratliff, L.~J.; and Sastry, S.~S. 2020.
\newblock On gradient-based learning in continuous games.
\newblock \emph{SIAM Journal on Mathematics of Data Science}, 2(1): 103--131.

\bibitem[{Mescheder, Nowozin, and Geiger(2017)}]{mescheder2017numerics}
Mescheder, L.; Nowozin, S.; and Geiger, A. 2017.
\newblock The numerics of gans.
\newblock \emph{arXiv preprint arXiv:1705.10461}.

\bibitem[{Nagarajan and Kolter(2017)}]{nagarajan2017gradient}
Nagarajan, V.; and Kolter, J.~Z. 2017.
\newblock Gradient descent GAN optimization is locally stable.
\newblock \emph{arXiv preprint arXiv:1706.04156}.

\bibitem[{Neumann(1928)}]{neumann1928theorie}
Neumann, J.~V. 1928.
\newblock Zur theorie der gesellschaftsspiele.
\newblock \emph{Mathematische annalen}, 100(1): 295--320.

\bibitem[{Nouiehed et~al.(2019)Nouiehed, Sanjabi, Huang, Lee, and
  Razaviyayn}]{nouiehed_solving_2019}
Nouiehed, M.; Sanjabi, M.; Huang, T.; Lee, J.~D.; and Razaviyayn, M. 2019.
\newblock Solving a Class of Non-Convex Min-Max Games Using Iterative First
  Order Methods.
\newblock In Wallach, H.; Larochelle, H.; Beygelzimer, A.; Alch{\'e}-Buc,
  F.~d.; Fox, E.; and Garnett, R., eds., \emph{Advances in Neural Information
  Processing Systems}, volume~32. Curran Associates, Inc.

\bibitem[{Polyak(1963)}]{polyak1963gradient}
Polyak, B.~T. 1963.
\newblock Gradient methods for minimizing functionals.
\newblock \emph{Zhurnal vychislitel'noi matematiki i matematicheskoi fiziki},
  3(4): 643--653.

\bibitem[{Robbins and Siegmund(1971)}]{robbins1971convergence}
Robbins, H.; and Siegmund, D. 1971.
\newblock A convergence theorem for non negative almost supermartingales and
  some applications.
\newblock In \emph{Optimizing methods in statistics}, 233--257. Elsevier.

\bibitem[{Rockafellar(1970)}]{rockafellar_convex_1970}
Rockafellar, R.~T. 1970.
\newblock \emph{Convex analysis}.
\newblock Princeton Mathematical Series, No. 28. Princeton University Press,
  Princeton, N.J.

\bibitem[{Sanjabi et~al.(2018)Sanjabi, Ba, Razaviyayn, and
  Lee}]{sanjabi2018convergence}
Sanjabi, M.; Ba, J.; Razaviyayn, M.; and Lee, J.~D. 2018.
\newblock On the convergence and robustness of training gans with regularized
  optimal transport.
\newblock \emph{arXiv preprint arXiv:1802.08249}.

\bibitem[{Sebbouh, Cuturi, and Peyr{\'e}(2022)}]{sebbouh2021randomized}
Sebbouh, O.; Cuturi, M.; and Peyr{\'e}, G. 2022.
\newblock Randomized Stochastic Gradient Descent Ascent.
\newblock In \emph{International Conference on Artificial Intelligence and
  Statistics}, 2941--2969. PMLR.

\bibitem[{Sharma et~al.(2022)Sharma, Panda, Joshi, and
  Varshney}]{sharma_2022_federated}
Sharma, P.; Panda, R.; Joshi, G.; and Varshney, P. 2022.
\newblock Federated Minimax Optimization: Improved Convergence Analyses and
  Algorithms.
\newblock In Chaudhuri, K.; Jegelka, S.; Song, L.; Szepesvari, C.; Niu, G.; and
  Sabato, S., eds., \emph{Proceedings of the 39th International Conference on
  Machine Learning}, volume 162 of \emph{Proceedings of Machine Learning
  Research}, 19683--19730. PMLR.

\bibitem[{Sinha et~al.(2017)Sinha, Namkoong, Volpi, and
  Duchi}]{sinha2017certifying}
Sinha, A.; Namkoong, H.; Volpi, R.; and Duchi, J. 2017.
\newblock Certifying some distributional robustness with principled adversarial
  training.
\newblock \emph{arXiv preprint arXiv:1710.10571}.

\bibitem[{Xian et~al.(2021)Xian, Huang, Zhang, and Huang}]{xian_2021_faster}
Xian, W.; Huang, F.; Zhang, Y.; and Huang, H. 2021.
\newblock A Faster Decentralized Algorithm for Nonconvex Minimax Problems.
\newblock In Ranzato, M.; Beygelzimer, A.; Dauphin, Y.; Liang, P.; and Vaughan,
  J.~W., eds., \emph{Advances in Neural Information Processing Systems},
  volume~34, 25865--25877. Curran Associates, Inc.

\bibitem[{Yan et~al.(2020)Yan, Xu, Lin, Liu, and Yang}]{yan2020optimal}
Yan, Y.; Xu, Y.; Lin, Q.; Liu, W.; and Yang, T. 2020.
\newblock Optimal epoch stochastic gradient descent ascent methods for min-max
  optimization.
\newblock \emph{arXiv preprint arXiv:2002.05309}.

\bibitem[{Yang et~al.(2021)Yang, Orvieto, Lucchi, and He}]{yang2021faster}
Yang, J.; Orvieto, A.; Lucchi, A.; and He, N. 2021.
\newblock Faster Single-loop Algorithms for Minimax Optimization without Strong
  Concavity.
\newblock \emph{arXiv preprint arXiv:2112.05604}.

\bibitem[{Zhang, Yang, and Ba{\c{s}}ar(2021)}]{zhang2021multi}
Zhang, K.; Yang, Z.; and Ba{\c{s}}ar, T. 2021.
\newblock Multi-agent reinforcement learning: A selective overview of theories
  and algorithms.
\newblock \emph{Handbook of Reinforcement Learning and Control}, 321--384.

\end{thebibliography}

\onecolumn
\begin{appendices}
\section{Proofs in Section \ref{sec: rsgda}}
\label{app: rsgda}

\subsection{Proofs in Section \ref{sec: 4.2}}

\begin{proof}[\textbf{Proof of Theorem \ref{thm: rgda is contractive}}]

For ease of notation, we define the operator $H$ as follows: 
\begin{equation}
    H (x, y) := 
    \begin{dcases}
        (x - \alpha \nabla_x F(x, y), y), & \quad \text{w.p.} \quad p, \\ 
        (x, y + \alpha \nabla_y F(x, y)), & \quad \text{w.p.} \quad 1 - p.
    \end{dcases}
\end{equation}
We also define 
\begin{equation}
    u = \begin{pmatrix}
    x \\ 
    y
    \end{pmatrix}, \quad 
    G(u) = \begin{pmatrix}
    \alpha \nabla_x F(x, y) \\
    - \alpha \nabla_y F(x, y)
    \end{pmatrix}.
\end{equation}
For simplicity, we denote $G(u_k)$ by $G_k$. 
Since $F$ is $\mu$-SCSC, we have 
\begin{equation}
    \frac{1}{\alpha} \langle G(u) - G(\tilde{u}), u - \tilde{u} \rangle \geq \mu \|  u - \tilde{u} \|^2, \quad \forall u, \tilde{u}.
\end{equation}
Particularly, if $u^*$ is the saddle point of $F$, we have
\begin{equation} \label{equ: a6}
    \langle G(u) , u - u^* \rangle \geq \alpha \mu \|  u - u^* \|^2, \quad \forall u.
\end{equation}

Recall that $\e_k[\cdot]$ is the conditional expectation with respect to the filter $\mathcal{F}_k = \{ x_0, x_1, \dots, x_k \}$, thus we have 
\begin{equation} \label{equ: a1}
\begin{aligned}
    \e_k \left[\| H (x_k, y_k) - (x^*, y^*) \|^2\right] &= p \| (x_k - \alpha \nabla_x F(x_k, y_k), y_k) - (x^*, y^*) \|^2 \\
    & \quad + (1 - p) \| (x_k, y_k + \alpha \nabla_y F(x_k, y_k)) - (x^*, y^*) \|^2.
\end{aligned}
\end{equation}
Now, we estimate the two terms separately. 
For the first term, we have
\begin{equation} \label{equ: a2}
\begin{aligned}
    \| (x_k - \alpha \nabla_x F(x_k, y_k), y_k) - (x^*, y^*) \|^2 & = \| x_k - \alpha \nabla_x F(x_k, y_k) - x^* \|^2 + \| y_k - y^* \|^2 \\ 
    & = \| x_k - x^* \|^2 + \| y_k - y^* \|^2 \\ & \quad +  \alpha^2 \| \nabla_x F(x_k, y_k) \|^2 - 2 \alpha \langle x_k - x^*, \nabla_x F(x_k, y_k) \rangle.
\end{aligned}
\end{equation}
Similarly, for the second term, we can obtain that 
\begin{equation} \label{equ: a3}
\begin{aligned}
     \| (x_k, y_k + \alpha \nabla_y F(x_k, y_k)) - (x^*, y^*) \|^2 &= \| x_k - x^* \|^2 + \| y_k - y^* \|^2 \\
     & \quad + \alpha^2 \| \nabla_y F(x_k, y_k) \|^2 + 2 \alpha \langle y_k - y^*, \nabla_y F(x_k, y_k) \rangle.
\end{aligned}
\end{equation}
Putting the above three equations, \eqref{equ: a1}, \eqref{equ: a2} and \eqref{equ: a3} together, we have 
\begin{equation} \label{equ: a4}
\begin{aligned}
    \e_k \left[\| H (x_k, y_k) - (x^*, y^*) \|^2\right] = & \| (x_k, y_k) - (x^*, y^*) \|^2 + \alpha^2 p \| \nabla_x F(x_k, y_k) \|^2 + \alpha^2 (1-p) \| \nabla_y F(x_k, y_k) \|^2 \\
    & - 2 \alpha \left( p \langle x_k - x^*, \nabla_x F(x_k, y_k) \rangle - (1-p) \langle y_k - y^*, \nabla_y F(x_k, y_k) \rangle \right).
\end{aligned}
\end{equation}
Without loss of generality, we assume that $p \leq 1 / 2$ \footnote{Note that $x$ and $y$ are symmetric in this algorithm, which means that if $p \geq 1/2$, one can exchange the notations $x$ and $y$ and define a new probability parameter $\tilde{p}:= 1 - p$ to complete the proof.} Now we analyze the equality~\eqref{equ: a4} under two cases. 
\begin{itemize}
    \item \underline{Case 1: $\langle y_k - y^*, \nabla_y F(x_k, y_k) \rangle \geq 0$.} By rearranging (\ref{equ: a4}), we have 
\begin{equation} 
\begin{aligned}
    \e_k \left[\| H (x_k, y_k) - (x^*, y^*) \|^2\right] = & \| u_k - u^* \|^2 + \alpha^2 (p \| \nabla_x F(x_k, y_k) \|^2 + (1-p) \| \nabla_y F(x_k, y_k) \|^2) \\
    & - 2 \alpha \left( p \langle x_k - x^*, \nabla_x F(x_k, y_k) \rangle - (1-p) \langle y_k - y^*, \nabla_y F(x_k, y_k) \rangle \right) \\ 
\end{aligned}
\end{equation}
Since $p \leq 1/2$, we can obtain that
\begin{equation} \label{eqn:ap}
\begin{aligned}
& \alpha (p \langle x_k - x^*, \nabla_x F(x_k, y_k) \rangle - (1-p) \langle y_k - y^*, \nabla_y F(x_k, y_k) \rangle) \\
& \qquad \qquad \ge p (\langle x_k - x^*, \alpha \nabla_x F(x_k, y_k) \rangle - \langle y_k - y^*, \alpha \nabla_y F(x_k, y_k) \rangle) \\
    & \qquad \qquad = p \langle u_k - u^*, G_k \rangle \ge p \alpha\mu \| u_k - u^* \|^2,
\end{aligned}
\end{equation}
where the last inequality is derived by~\eqref{equ: a6},
and 
\begin{equation}
\begin{aligned}
    & \alpha^2(p \| \nabla_x F(x_k, y_k) \|^2 + (1-p) \| \nabla_y F(x_k, y_k) \|^2) \\
    & \qquad \qquad \leq (1-p) \left(\| \alpha \nabla_x F(x_k, y_k) \|^2 + \| \alpha \nabla_y F(x_k, y_k)\| \right) \\
    & \qquad \qquad  =  (1-p) \|G_k\|^2 \leq (1-p) \alpha^2 L_1 \|u_k - u^*\|^2, 
\end{aligned}
\end{equation}
where the last inequality is derived by the $L_1$ smoothness of $F$. 
Therefore, we can obtain that
\begin{equation}
% \begin{aligned}
    \e_k \left[\| H (x_k, y_k) - (x^*, y^*) \|^2 \right] 
    % & \| u_k - u^* \|^2  - 2 p \mu \alpha \| u_k - u^* \|^2 
    % & + \alpha^2 (1-p) \| \nabla_x F(x_k, y_k) \|^2 + \alpha^2 (1-p) \| \nabla_y F(x_k, y_k) \|^2 \\
    \leq (1 + \alpha^2 (1 - p) L_1^2 - 2p \mu \alpha) \| u_k - u^* \|^2.
% \end{aligned}
\end{equation}
Then, for $\alpha$ sufficiently small (i.e., $0 < \alpha  < \frac{2p\mu}{(1-p)L_1^2}$), we can define $\rho := 1 - 2p \mu \alpha + \alpha^2 (1 - p) L_1^2 < 1$, and thus 
\begin{equation}
    \e_k \left[\| H (x_k, y_k) - (x^*, y^*) \|^2\right] \leq \rho \| (x_k, y_k) - (x^*, y^*) \|^2.
\end{equation}

\item \underline{Case 2: $\langle y_k - y^*, \nabla_y F(x_k, y_k) \rangle \leq 0$.} 
In this case, the inequality~\eqref{eqn:ap} does not hold any more. 
To handle this issue, we decompose $(1-p)$ as $p+(1-2p)$, and obtain that
\begin{equation}
\begin{aligned}
& \alpha (p \langle x_k - x^*, \nabla_x F(x_k, y_k) \rangle - (1-p) \langle y_k - y^*, \nabla_y F(x_k, y_k) \rangle) \\
& \qquad \qquad \ge p (\langle x_k - x^*, \alpha \nabla_x F(x_k, y_k) \rangle - \langle y_k - y^*, \alpha \nabla_y F(x_k, y_k) \rangle) - (1-2p)\| \langle y_k - y^*, \alpha \nabla_y F(x_k, y_k) \rangle \| \\ 
& \qquad \qquad \ge p \alpha\mu \| u_k - u^* \|^2 - (1-2p) \langle y_k - y^*, \alpha \nabla_y F(x_k, y_k) \rangle \ge p \alpha\mu \| u_k - u^* \|^2, 
\end{aligned}
\end{equation}
where the last inequality is due to that $F$ is strongly concave of $y$. 
% According to~\eqref{equ: a6}, we have 
% \begin{equation} \label{equ: a9}
    % \langle x_k - x^*, \nabla_x F(x_k, y_k) \rangle \geq \mu \| u_k - u^* \|^2 + \langle y_k - y^*, \nabla_y F(x_k, y_k) \rangle.
% \end{equation}
Therefore, the inequality 
\begin{equation}
% \begin{aligned}
    \e_k \left[\| H (x_k, y_k) - (x^*, y^*) \|^2\right] 
    % = & \| u_k - u^* \|^2 + \alpha^2 p \| \nabla_x F(x_k, y_k) \|^2 + \alpha^2 (1-p) \| \nabla_y F(x_k, y_k) \|^2 
    % \leq & \| u_k - u^* \|^2 -2 \alpha p \mu \| u_k - u^* \|^2 \\ 
    % & + \alpha^2 p \| \nabla_x F(x_k, y_k) \|^2 + \alpha^2 (1-p) \| \nabla_y F(x_k, y_k) \|^2 \\ 
    \leq (1 - 2 \alpha p \mu + \alpha^2 (1-p) L_1^2) \| u_k - u^* \|^2.
% \end{aligned}
\end{equation}
still holds for this case. 
Hence, for $\alpha$ sufficiently small, we also have
\begin{equation}
    \e_k \left[\| H (x_k, y_k) - (x^*, y^*) \|^2\right] \leq \rho \| (x_k, y_k) - (x^*, y^*) \|^2.
\end{equation}
\end{itemize} 

Combining the Case 1 and 2, we can conclude that,
for sufficiently small $\alpha$ ($0 < \alpha  < \frac{2p\mu}{(1-p)L_1^2}$), and $\rho := (1 - 2 \alpha \mu + \alpha^2 (1 - p) L_1^2 < 1$, 
we have 
\begin{equation}
    \e_k \left[\| H (x_k, y_k) - (x^*, y^*) \|^2\right] \leq \rho \| (x_k, y_k) - (x^*, y^*) \|^2.
\end{equation}

\end{proof}

A direct observation of Theorem \ref{thm: rgda is contractive} indicates the Corollary \ref{coro: 1}. Precisely, we provide the following proof. 

\begin{proof}[\textbf{Proof of Corollary \ref{coro: 1}}]

Let $\{ (x_k, y_k) \}$ be the sequence generated by Algorithm \ref{algo: rsgda}. By Theorem \ref{thm: rgda is contractive}, for $\alpha$ sufficiently small, there exists constant $\rho < 1$, such that  
\begin{equation} \label{equ: a12}
    \e_k [\| H (x_k, y_k) - (x^*, y^*) \|^2] \leq \rho \| (x_k, y_k) - (x^*, y^*) \|^2, \quad \forall k .
\end{equation}
Taking expectation in both sides of (\ref{equ: a12}), we have 
\begin{equation} 
    \e \left[\| H (x_k, y_k) - (x^*, y^*) \|^2 \right] \leq \rho \e \| (x_k, y_k) - (x^*, y^*) \|^2, \quad \forall k \geq 0.
\end{equation}
Hence that 
\begin{equation} 
    \e \left[\| H (x_k, y_k) - (x^*, y^*) \|^2 \right] \leq \rho^k \| (x_0, y_0) - (x^*, y^*) \|^2, \quad \forall k \geq 0,
\end{equation}
which indicates that $\{ (x_k, y_k) \}_k$ linearly converges to the saddle point $(x^*, y^*)$ in expectation.

\end{proof}

\subsection{Proofs in Section \ref{sec: 4.3}}

To prove the results in Section \ref{sec: 4.3}, we introduce three technical lemmas first.

\begin{lemma}[\citet{nouiehed_solving_2019}] \label{lemma: 1}
Under  Assumption \ref{assum: smoothness} and \ref{assum: ncpl condtion}, $\phi$ is $L_2$ smooth with $L_2 = L_1 + \frac{L_1 \kappa}{2}$ and $\kappa = L_1 / \mu$.
\end{lemma}

\begin{lemma}[\citet{karimi2016linear}] \label{lemma: 2}
    
If $g(\cdot)$ is $l$-smooth and it satisfies P{\L} condition with constant $\mu$, i.e., 
\begin{equation*}
    \| \nabla g(x) \|^2 \geq 2 \mu [g(x) - \min_x g(x)], \quad \forall x, 
\end{equation*}
then it also satisfies error bound condition with $\mu$, i.e., 
\begin{equation*}
    \| \nabla g(x) \|^2 \geq \mu \| x_p - x \|, \quad \forall x,
\end{equation*}
where $x_p$ stands for the projection of $x$ onto the optimal set.
\end{lemma}

\begin{lemma}[Robbins-Siegmund theorem \citep{robbins1971convergence}] \label{lemma: 3}
    
Consider a filtration $\{ \mathcal{F}_k \}_k$, the nonnegative sequences of $\{ \mathcal{F}_k \}_k$ adapted processes $\{ V_k \}_k$, $\{ U_k \}_k$ and $\{ Z_k \}_k$ such that $\sum_k Z_k < +\infty$ almost surly, and 
\begin{equation}
    \e_k [V_{k+1} |\mathcal{F}_k ] + U_{k+1} \leq V_k + Z_k, \quad \forall k \geq 0.
\end{equation}
Then $\{ V_k \}_k$ converges and $\sum_k U_k < +\infty$ almost surly.
    
\end{lemma}

The basic idea to prove Theorem \ref{thm: convergence of rsgda} is to introduce a suitable decreasing Lyapunov function $V$, such that $\{ V (x_k, y_k) \}_k$ is decreasing. In this paper, we introduce the following function. 
\begin{equation} \label{equ: a15}
    V_k := V(x_k, y_k) = \phi (x_k) + C[\phi (x_k) - F(x_k, y_k)] = (1+C) \phi (x_k) - C F(x_k, y_k),
\end{equation}
where $C > 0$ is a constant to be determined later. 

\begin{proof}[\textbf{Proof of Theorem \ref{thm: convergence of rsgda}}]

Let $V_k$ be the function given by (\ref{equ: a15}). Without loss of generality, we assume that $\phi (x_k) \geq 0$, $\forall k \geq 0$. Then $\{ V_k \}_k$ is a nonnegative sequence. Next we need to estimate the value 
\begin{equation}
    V_k - \e_k \left[V_{k+1}\right].
\end{equation}
Here we denote $\e_{z_k}[\cdot]$ by the expectation conditioned on the random variable $z_k$, and $\e_k[\cdot]$ by the expectation conditioned on all past random variables.

\begin{itemize}
    \item \underline{Step 1: Estimate $\e_k \left[\phi (x_{k+1})\right] - \phi (x_k)$.} Recall the definition of $x_{k+1}$, we have 
\begin{equation} \label{equ: a17}
    \e_k \left[\phi (x_{k+1})\right] = p \e_{z_k} \left[\phi (x_k^+)\right] + (1-p) \e_{z_k}\left[\phi (x_k)\right] = p \e_{z_k} \left[\phi (x_k^+)\right] + (1-p) \phi (x_k).
\end{equation}
Based on the Assumption \ref{assum: smoothness} and \ref{assum: ncpl condtion} hold, then $\phi$ is $L_2$-smooth by Lemma \ref{lemma: 1}. Hence, 
\begin{equation}
\begin{aligned}
    \phi (x_k^+) & \leq \phi (x_k) + \langle x_k^+ - x_k, \nabla \phi (x_k) \rangle + \frac{L_2}{2} \| x_k^+ - x_k \|^2 \\ 
    & = \phi (x_k) - \alpha_k \langle \nabla_x f(x_k, y_k; z_k), \nabla \phi (x_k) \rangle + \frac{\alpha_k^2 L_2}{2} \| \nabla_x f(x_k, y_k; z_k) \|^2.
\end{aligned}
\end{equation}
Taking conditional expectation in both sides of the above equation, we obtain
\begin{equation} \label{equ: a18}
\begin{aligned}
    \e_{z_k}\left[\phi (x_k^+)\right] & \leq \phi (x_k) - \alpha_k \langle \nabla_x F(x_k, y_k), \nabla \phi (x_k) \rangle + \frac{\alpha_k^2 L_2}{2} \e_{z_k} \left[\| \nabla_x f(x_k, y_k; z_k) \|^2\right] \\
    & = \phi (x_k) - \alpha_k \langle \nabla_x F(x_k, y_k), \nabla \phi (x_k) \rangle + \frac{\alpha_k^2 L_2}{2} \e_{z_k} \left[\| \nabla_x F(x_k, y_k) - \nabla f(x_k, y_k; z_k) \|^2 \right] \\ 
    & \quad + \frac{\alpha_k^2 L_2}{2} \e_{z_k} \left[\| \nabla_x F(x_k, y_k) \|^2\right] + \alpha_k^2 L_2 \e_{z_k} \left[\langle \nabla_x F(x_k, y_k), \nabla_x F(x_k, y_k) - \nabla f(x_k, y_k; z_k) \right]. 
\end{aligned}
\end{equation}
By using Assumption 2, we can bound and eliminate the last two terms in the equation, respectively. Thus, 
\begin{equation}
\e_{z_k}\left[\phi (x_k^+)\right] \leq  \phi (x_k) - \alpha_k \langle \nabla_x F(x_k, y_k), \nabla \phi (x_k) \rangle + \frac{\alpha_k^2 L_2 \sigma^2}{2} + \frac{\alpha_k^2 L_2}{2} \| \nabla_x F(x_k, y_k) \|^2. 
\end{equation}
Furthermore, by rewritting the term $\langle \nabla_x F(x_k, y_k), \nabla \phi (x_k) \rangle$, we have
\begin{equation} \label{equ:res}
\begin{aligned}
\e_{z_k}\left[\phi (x_k^+)\right]
     & \leq \phi (x_k)  + \frac{\alpha_k^2 L_2 \sigma^2}{2} + \frac{\alpha_k}{4} \| \nabla_x F(x_k, y_k) \|^2 \\ 
     & \quad + \frac{\alpha_k}{2} \| \nabla_x F(x_k, y_k) - \nabla \phi (x_k) \|^2 - \frac{\alpha_k}{2} \| \nabla \phi (x_k) \|^2 - \frac{\alpha_k}{2} \| \nabla_x F(x_k, y_k) \|^2 \\
    & =  \phi (x_k) - \frac{\alpha_k}{2} \| \nabla \phi (x_k) \|^2 - \frac{\alpha_k}{4} \| \nabla_x F(x_k, y_k) \|^2 + \frac{\alpha_k}{2} \| \nabla_x F(x_k, y_k) - \nabla \phi (x_k) \|^2 + \frac{\alpha_k^2 L_2 \sigma^2}{2}.
\end{aligned}
\end{equation}
Combining equation~\eqref{equ: a17} and  ~\eqref{equ:res}, we have 
\begin{equation}
    \e_k \left[\phi (x_{k+1})\right] \leq \phi (x_k) - \frac{\alpha_k p}{2} \| \nabla \phi (x_k) \|^2 - \frac{\alpha_k p}{4} \| \nabla_x F(x_k, y_k) \|^2 + \frac{\alpha_k p}{2} \| \nabla_x F(x_k, y_k) - \nabla \phi (x_k) \|^2 + \frac{\alpha_k^2 L_2 \sigma^2 p}{2}.
\end{equation}

\item \underline{Step 2: Estimate the term $\e_k [F(x_{k+1}, y_{k+1})] - F(x_k, y_k)$.} Note that 
\begin{equation} \label{equ: a20}
    \e_k [F(x_{k+1}, y_{k+1})] - F(x_k, y_k) = p \left(\e_{z_k} [F(x_k^+, y_k)] - F(x_k, y_k) \right) + p \left(\e_{z_k} [F(x_k, y_k^+)] - F(x_k, y_k) \right).
\end{equation}
Now we estimate the two terms separately. First, because $F$ is $L_1$-smooth, we obtain 
\begin{equation*}
\begin{aligned}
    F(x_k^+, y_k) & \geq F(x_k, y_k) + \langle \nabla_x F(x_k, y_k), x_k^+ - x_k \rangle - \frac{L_1}{2} \|  x_k^+ - x_k \|^2 \\ 
    & = F(x_k, y_k) - \alpha_k \langle \nabla_x F(x_k, y_k), \nabla_x f(x_k, y_k; z_k) \rangle - \frac{\alpha_k^2 L_1}{2} \|   \nabla_x f(x_k, y_k; z_k) \|^2.
\end{aligned}
\end{equation*}
Hence, based on Assumption 2, we can obtain that
\begin{equation} \label{equ: a21}
\begin{aligned}
    \e_{z_k} \left[ F(x_k^+, y_k) \right] - F(x_k, y_k) 
    & \geq - \alpha_k \| \nabla_x F(x_k, y_k) \|^2 - \frac{\alpha_k^2 L_1}{2} \| \nabla_x F(x_k, y_k) \|^2 \\
    & \quad - \frac{\alpha_k^2 L_1}{2} \e_{z_k} \left[ \| \nabla_x f(x_k, y_k; z_k) - \nabla_x F(x_k, y_k) \|^2 \right] .
\end{aligned}
\end{equation}
Furthermore, since $\alpha_k \leq 1/L_2 < 1/L_1$, we have
\begin{equation}
\e_{z_k} \left[ F(x_k^+, y_k) \right] - F(x_k, y_k) \ge - \frac{3 \alpha_k}{2} \| \nabla_x F(x_k, y_k) \|^2 - \frac{\alpha_k^2 L_1 \sigma^2}{2}.
\end{equation}
Similarly, we have
\begin{equation*}
\begin{aligned}
    F(x_k, y_k^+) & \geq F(x_k, y_k) + \langle \nabla_y F(x_k, y_k), y_k^+ - y_k \rangle - \frac{L_1}{2} \|  y_k^+ - y_k \|^2 \\ 
    & = F(x_k, y_k) + \eta_k \langle \nabla_y F(x_k, y_k), \nabla_y f(x_k, y_k; z_k) \rangle - \frac{\eta_k^2 L_1}{2} \|   \nabla_y f(x_k, y_k; z_k) \|^2.
\end{aligned}
\end{equation*}
Hence, 
\begin{equation} \label{equ: a22}
\begin{aligned}
    \e_{z_k} \left[ F(x_k, y_k^+) \right] - F(x_k, y_k) & \geq \eta_k \| \nabla_y F(x_k, y_k) \|^2 - \frac{\eta_k^2 L_1}{2} \| \nabla_y F(x_k, y_k) \|^2 - \frac{\eta_k^2 L_1 \sigma^2}{2} \\ 
    & \geq \frac{\eta_k}{2} \| \nabla_y F(x_k, y_k) \|^2 - \frac{\eta_k^2 L_1 \sigma^2}{2}.
\end{aligned}
\end{equation}
Combining~\eqref{equ: a20}, \eqref{equ: a21} and \eqref{equ: a22}, we obtain 
{\small
\begin{equation}
\begin{aligned}
    \e_k [F(x_{k+1}, y_{k+1})] - F(x_k, y_k) & \geq p \left( - \frac{3\alpha_k}{2} \| \nabla_x F(x_k, y_k) \|^2 - \frac{\alpha_k^2 L_1 \sigma^2}{2} \right)  + (1-p) \left( \frac{\eta_k}{2} \| \nabla_y F(x_k, y_k) \|^2 - \frac{\eta_k^2 L_1 \sigma^2}{2} \right) \\ 
    & = \frac{(1-p)\eta_k}{2} \| \nabla_y F(x_k, y_k) \|^2 -\frac{3\alpha_k p}{2} \| \nabla_x F(x_k, y_k) \|^2 - \frac{\sigma^2L_1 (p \alpha_k^2 + (1-p) \eta_k^2)}{2}
\end{aligned}
\end{equation}
}
\item \underline{Step 3: Upper bound of $\| \nabla_x F(x_k, y_k) - \phi (x) \|^2$.}
Recall $y^* (x) := \arg\max_y F(x, y)$, thus $\nabla_y F(x, y^* (x)) = 0$ can be derived by Fermat's rule. Now, the gradient $\nabla \phi(x)$ can be computed as
\begin{equation*}
    \nabla \phi (x) = \nabla_x F(x, y^*(x)) = \nabla_x F(x, y^*(x)) + \nabla_y F(x, y^*(x)) \nabla_x y^*(x) = \nabla_x F(x, y^*(x)).
\end{equation*}
Furthermore, based on the smoothness of $F$ and the P{\L} condition, we have 
\begin{equation} \label{equ: a23}
\begin{aligned}
    \| \nabla_x F(x_k, y_k) - \nabla \phi (x_k) \|^2 & = \| \nabla_x F(x_k, y_k) - \nabla_x F(x_k, y^*(x_k)) \|^2 \\
    & \leq L_1^2 \| y_k - y^*(x_k) \|^2 
    \leq \left( \frac{L_1}{\mu} \right)^2 \| \nabla_y F(x_k, y_k) \|^2,
\end{aligned}
\end{equation}
where the last inequality holds due to the Lemma \ref{lemma: 2}.

\item \underline{Step 4: Estimate the Lyapunov function $V_k$.}
First, note that 
\begin{equation}
    \| a \|^2 + \| a - b \|^2 \geq \frac{1}{2} \| b \|^2, \quad \forall a, b
\end{equation}
according to the Young's inequality. Hence, 
\begin{equation*}
    -\frac{3p \alpha_k}{2} \| \nabla_x F(x_k, y_k) \|^2 \geq -3 p \alpha_k \| \nabla \phi (x_k) \|^2 - 3 p \alpha_k \| \nabla \phi (x_k) - \nabla_x F(x_k, y_k) \|^2.
\end{equation*}
Combining Step 1 and 2, we have 
\begin{equation} \label{equ: a24}
\begin{aligned}
V_{k}-\mathbb{E}_{k}\left[V_{k+1}\right]
& =(1+C)\left(\phi\left(x_{k}\right)-\mathbb{E}_{k}\left[\phi\left(x_{k+1}\right)\right]\right)+C\left(\mathbb{E}_{k}\left[F\left(x_{k+1}, y_{k+1}\right)\right]-F\left(x_{k}, y_{k}\right)\right) \\
& \geq (1+C)\left(\frac{p \alpha_{k}}{2}\left\|\nabla \phi\left(x_{k}\right)\right\|^{2}-\frac{p \alpha_{k}}{2}\left\|\nabla \phi\left(x_{k}\right)-\nabla_{x} F\left(x_{k}, y_{k}\right)\right\|^{2}+\frac{p \alpha_{k}}{4}\left\|\nabla_{x} F\left(x_{k}, y_{k}\right)\right\|^{2}\right) \\
& \qquad + C\left(\frac{(1-p) \eta_{k}}{2}\left\|\nabla_{y} F\left(x_{k}, y_{k}\right)\right\|^{2}-\frac{3 p \alpha_{k}}{2}\left\|\nabla_{x} F\left(x_{k}, y_{k}\right)\right\|^{2}\right) + U_k (C) \\
& \geq  {\left(\frac{p \alpha_{k}}{2}(1+C)-3 p \alpha_{k} C\right) \cdot\left\|\nabla \phi\left(x_{k}\right)\right\|^{2} -\left(\frac{p \alpha_{k}}{2}(1+C)+3 p \alpha_{k} \cdot C\right) \cdot\left\|\nabla_{x} F\left(x_{k}, y_{k}\right)-\nabla \phi \left(x_{k}\right)\right\|^{2} } \\ 
& \qquad +\frac{p \alpha_{k}}{4}(1+C)\left\|\nabla_{x} F\left(x_{k}, y_{k}\right)\right\|^{2}+\frac{(1-p) \eta_{k}}{2} C \cdot\left\|\nabla_{y} F\left(x_{k}, y_{k}\right)\right\|^{2} + U_k (C).
% - \left[\frac{L_{2} \alpha_{k}^{2} p(1+C)}{2}+\frac{p}{2} L_{1} \alpha_{k}^{2}+\frac{1-p}{2} L_{1} \eta_{k}^{2}\right] \sigma^{2} . 
\end{aligned}
\end{equation}
where 
\begin{equation}
    U_k(C) = -\left(\frac{L_{2} \alpha_{k}^{2} p(1+C)}{2}+\frac{p}{2} L_{1} \alpha_{k}^{2}+\frac{1-p}{2} L_{1} \eta_{k}^{2}\right) \sigma^{2}.
\end{equation}
Plugging~\eqref{equ: a23} into~\eqref{equ: a24}, we get
\begin{equation}
\begin{aligned}
V_{k}-\mathbb{E}_{k}\left[V_{k+1}\right] 
& \ge {\left(\frac{p \alpha_{k}}{2}(1+C)-3 p \alpha_{k} C\right) \cdot\left\|\nabla \phi\left(x_{k}\right)\right\|^{2}-\left(\frac{p \alpha_{k}}{2}(1+C)+3 p \alpha_{k} \cdot C\right)\left(\frac{L_{1}}{\mu}\right)^{2} \cdot\left\|\nabla_{y} F\left(x_{k}, y_{k}\right)\right\|^{2} } \\
& \qquad + \frac{p \alpha_{k}}{4}(1+C)\left\|\nabla_{x} F\left(x_{k}, y_{k}\right)\right\|^{2}+\frac{(1-p) \eta_{k}}{2} C \cdot\left\|\nabla_{y} F\left(x_{k}, y_{k}\right)\right\|^{2} + U_k (C) \\
% &-\left[\frac{L_{2} \alpha_{k}^{2} p(1+C)}{2}+\frac{p}{2} L_{1} \alpha_{k}^{2}+\frac{1-p}{2} L_{1} \eta_{k}^{2}\right] \sigma^{2} . \\
& = \left(\frac{(1-p) \eta_{k}}{2} C-\left(\frac{p \alpha_{k}}{2}(1+C)+3 p \alpha_{k} C\right) \cdot\left(\frac{L_{1}}{\mu}\right)^{2}\right) \cdot\left\|\nabla_{y} F\left(x_{k}, y_{k}\right)\right\|^{2} \\
& \qquad + {\left(\frac{p \alpha_{k}}{2}(1+C)-3 p \alpha_{k} C\right) \cdot \| \nabla \phi (x_{k}) \|^{2} } +\frac{p \alpha_{k}}{4}(1+C)\left\|\nabla_{x} F\left(x_{k}, y_{k}\right)\right\|^{2} + U_k (C).
\end{aligned}
\end{equation}
Now, we try to determine the concrete value of $C$. 
Here, we choose $C = 1/ 10$, then 
\begin{equation*}
\begin{aligned}
& \frac{p \alpha_{k}}{2}(1+C)-3 p \alpha_{k} C = p \alpha_k \frac{1-5C}{2} = \frac{1}{4} p \alpha_k, \\ 
& \frac{(1-p) \eta_{k}}{2} C-\left(\frac{p \alpha_{k}}{2}(1+C)+3 p \alpha_{k} C\right) \cdot\left(\frac{L_{1}}{\mu}\right)^{2} \geq \frac{p}{20} \left(\frac{L_{1}}{\mu}\right)^{2} \alpha_k.
\end{aligned}
\end{equation*}
Hence, 
\begin{equation} \label{equ: a26}
\begin{aligned}
V_{k}-\mathbb{E}_{k}\left[V_{k+1}\right] \ge \frac{1}{4} p \alpha_{k} \cdot\left\|\nabla \phi\left(x_{k}\right)\right\|^{2}+\frac{p}{20}\left(\frac{L_{1}}{\mu}\right)^{2} \alpha_{k} \left\|\nabla_{y} F\left(x_{k}, y_{k}\right)\right\|^{2}+\frac{11}{40} p \alpha_{k} \left\|\nabla_{x} F\left(x_{k}, y_{k}\right)\right\|^{2} + U_k (\frac{1}{10})
% & -\left[\frac{11}{20} L_{2} p \alpha_{k}^{2}+\frac{p}{2} L_{1} \alpha_{k}^{2}+\frac{1-p}{2} L_{1} \eta_{k}^{2}\right] \cdot \sigma^{2} .
\end{aligned}
\end{equation}
Rearranging equation~\eqref{equ: a26}, we have 
\begin{equation} \label{equ: a27}
\begin{aligned}
\frac{1}{4} p \alpha_{k} \left\|\nabla \phi\left(x_{k}\right)\right\|^{2}+\frac{p}{20}\left(\frac{L_{1}}{\mu}\right)^{2} \alpha_{k} \left\|\nabla_{y} F\left(x_{k}, y_{k}\right)\right\|^{2}+\frac{11}{40} p\alpha_{k} \left\|\nabla_{x} F\left(x_{k}, y_{k}\right)\right\|^{2}+\mathbb{E}_{k}\left[V_{k+1}\right] 
\le & V_{k} + U_k (\frac{1}{10})
% +\left[\frac{11}{20} L_{2} p \alpha_{k}^{2}+\frac{p}{2} L_{1} \alpha_{k}^{2}+\frac{1-p}{2} L_{1} \eta_{k}^{2}\right] \cdot \sigma^{2}
\end{aligned}
\end{equation}
\item \underline{Step 5: Using Robbins-Siegmumd Theorem.}
In this step, we use a similar proof technique to \citet{sebbouh2021randomized}. 
We first define that, for all $k \geq 0$,
\begin{equation}
\begin{aligned}
& w_{k}:= \frac{2 \alpha_{k}}{\sum_{j=0}^{k} \alpha_{j}}, \quad 
g_{0}=\frac{1}{4}\left\|\nabla \phi\left(x_{0}\right)\right\|^{2}+\frac{1}{20}\left(\frac{L_{1}}{\mu}\right)^{2} \cdot\left\|\nabla_{y} F\left(x_{0}, y_{0}\right)\right\|^{2}+\frac{11}{40}\left\|\nabla_{x} F\left(x_{0}, y_{0}\right)\right\|^{2} , \\
& g_{k+1}:=\left(1-w_{k}\right) g_{k}+w_{k} \cdot\left(\frac{1}{4}\left\|\nabla \phi\left(x_{k}\right)\right\|^{2}+\frac{1}{20}\left(\frac{L_{1}}{\mu}\right)^{2}\left\|\nabla_{y} F\left(x_{k}, y_{k}\right)\right\|^{2}+\frac{11}{40}\left\|\nabla_{x} F\left(x_{k}, y_{k}\right)\right\|^{2}\right) .
\end{aligned}
\end{equation}
Since $\{ \alpha_k \}$ is non-increasing, we have $w_k \in [0, 1]$ for all $k \ge 0$. 
Hence, $g_k$ can be viewed as a convex combination of $\{ h_0, h_1, \dots, h_{k-1} \}$.

Next, using Lemma \ref{lemma: 3} for~\eqref{equ: a27}, and the fact that $\sum_k \alpha_k^2 < \infty, \sum_k \eta_k^2 < \infty$, we can conclude that $\{ V_k \}_k$ converges almost surly. 
On the other side, Rearranging~\eqref{equ: a27}, we have 
\begin{equation}
\frac{\sum_{j=0}^{k} \alpha_{j}}{2} g_{k+1}+\frac{1}{p} \mathbb{E}_{k}\left[V_{k+1}\right]+\frac{\alpha_{k}}{2} g_{k} \leqslant \frac{\sum_{j=0}^{k-1} \alpha_{j}}{2} g_{k}+\frac{1}{p} V_{k} +4 \cdot  U_k (\frac{1}{10}).
% \left[\frac{11}{20} L_{2} \alpha_{k}^{2}+\frac{1}{2} L_{1} \alpha_{k}^{2}+\frac{1-p}{2 p} L_{1} \eta_{k}^{2}\right] \sigma^{2} .
\end{equation}
Using Lemma \ref{lemma: 3} again and the fact that $\{ V_k \}_k$ converges almost surly, it can be obtained that $\{ \sum_{j=0}^{k} \alpha_{j} g_{k+1} \}_k$ converges almost surly and $\sum_k \alpha_k g_k < + \infty$. 
Particularly, this implies that $\lim_k \alpha_k g_k = 0$. 
Notice that $\alpha_k g_k = \frac{\alpha_k}{\sum_{j=0}^{k-1} \alpha_j} \sum_{j=0}^{k-1} \alpha_j g_k$, together with this $\{ \sum_{j=0}^{k} \alpha_{j} g_{k+1} \}_k$ converges almost surly and $\sum_k \frac{\alpha_k}{\sum_{j=0}^{k-1} \alpha_j} = \infty$, we have $\lim_k \sum_{j=0}^{k-1} g_k = 0$, i.e., 
\begin{equation*}
    g_k = o\left( \frac{1}{\sum_{j=0}^{k-1} \alpha_j} \right).
\end{equation*}
Finally, since $g_k$ is a convex combination of $\{ h_0, h_1, \dots, h_{k-1} \}$, we have that 
\begin{equation}
    \min_{t=0, 1, \dots, k-1} h_t \leq g_k = o\left( \frac{1}{\sum_{j=0}^{k-1} \alpha_j} \right).
\end{equation}
\end{itemize}
\end{proof}

Theorem \ref{thm: convergence rate} is a direct consequence of Theorem \ref{thm: convergence of rsgda}. The proof is Theorem \ref{thm: convergence rate} relies on the convergence rate of $1 / \sum_{j=0}^k \alpha_j$.

\begin{proof}[\textbf{Proof of Theorem \ref{thm: convergence rate}}]

According to Theorem \ref{thm: convergence of rsgda}, we have 
\begin{equation}
    \min_{t=0, \dots, k} h_t = o \left( \frac{1}{\sum_{j=0}^k \alpha_j} \right),
\end{equation}
almost surly. 
On the other side, for any $\epsilon > 0$ given, we define $\alpha_0 = 1$ and $\alpha_j = j^{-\frac{1}{2} - \epsilon}$ for $j = 1, 2, \dots$. Hence, $\sum_{j=0}^k \alpha_j = O(k^{\frac{1}{2} - \epsilon})$. Thus, we get that 
\begin{equation}
    \min_{t=0, \dots, k} h_t = o \left( \frac{1}{\sum_{j=0}^k \alpha_j} \right) = o (k^{- \frac{1}{2} + \epsilon}).
\end{equation}

\end{proof}

The proof of Corollary \ref{coro: 2} can be derived directly. 

\begin{proof}[\textbf{Proof of Corollary \ref{coro: 2}}]

Without loss of generality, we assume that $V_k \geq 0$ for all $k \geq 0$. The convergence in expectation proof is based on the telescopic cancellation in \eqref{equ: a27}.  Taking the expectation to \eqref{equ: a27} and summing for $k=0,1, \dots, n$, we obtain that 
\begin{equation}
    (n+1) \alpha  \cdot \min_{k=0,\dots, n} \e [h_k] \leq \alpha \sum_{k=0}^n \e [h_k] \leq \e [V_0] - \e [V_{n+1}] \leq \e [V_0], \quad \forall n \geq 0.
\end{equation}
Hence, 
\begin{equation}
    \min_{k=0,\dots, n} \e [h_k] \leq \frac{\e[V_0]}{(n+1) \alpha} = O \left(\frac{1}{n+1} \right).
\end{equation}

\end{proof}

Finally, we prove Corollary \ref{coro: 3} in the following.

\begin{proof}[\textbf{Proof of Corollary \ref{coro: 3}}]

The proof of Corollary \ref{coro: 3} is similar to the proof of Theorem \ref{thm: convergence of rsgda}. Without loss of generality, we assume that $V_k \geq 0$ for all $k \geq 0$. Following the same steps as the proof of Theorem \ref{thm: convergence of rsgda}, we can obtain that \eqref{equ: a27} holds for any $k \geq 0$, i.e., 
\begin{equation} \label{equ: a42}
\begin{aligned}
\frac{1}{4} p \alpha_{k} \left\|\nabla \phi\left(x_{k}\right)\right\|^{2}+\frac{p}{20}\left(\frac{L_{1}}{\mu}\right)^{2} \alpha_{k} \left\|\nabla_{y} F\left(x_{k}, y_{k}\right)\right\|^{2}+\frac{11}{40} p\alpha_{k} \left\|\nabla_{x} F\left(x_{k}, y_{k}\right)\right\|^{2}+\mathbb{E}_{k}\left[V_{k+1}\right] 
\le V_{k} + U_k (\frac{1}{10}).
% +\left[\frac{11}{20} L_{2} p \alpha_{k}^{2}+\frac{p}{2} L_{1} \alpha_{k}^{2}+\frac{1-p}{2} L_{1} \eta_{k}^{2}\right] \cdot \sigma^{2}, \quad \forall k \geq 0.
\end{aligned}
\end{equation}
Taking expectation to equation\eqref{equ: a42} and summing for $k=0,1,2, \dots, n$, we get that 
\begin{equation} \label{equ: a43}
    \alpha \sum_{k=0}^n \e [h_k] \leq \e [V_0] - \e [V_{n+1}] + (n+1) \left[\frac{11}{20} L_{2} p \alpha^{2}+\frac{p}{2} L_{1} \alpha^{2}+\frac{1-p}{2} L_{1} \frac{18^2 p^2}{(1-p)^2} \kappa^4 \alpha^2 \right] \cdot \sigma^{2}.
\end{equation}
For ease of notation, we define $M:= \frac{11}{20}L_2 p + \frac{p}{2}L_1 + \frac{18^2(1-p)p^2}{2(1-p)^2}L_1 \kappa^4$.
Hence, 
\begin{equation}
\begin{aligned}
\min_{k=0,\dots, n} \e [h_k] \leq \frac{1}{n+1} \sum_{k=0}^n \e [h_k]  \leq \frac{\e [V_0]}{(n+1) \alpha} + M \sigma^2 \alpha  = \frac{\sqrt{M \sigma^2 \e [V_0]}}{\sqrt{n+1}},
\end{aligned}
\end{equation}
where $\alpha = \sqrt{\frac{\e [V_0]}{(n+1)M}} \sigma$. Thus, for any $\epsilon > 0$, if $n \geq \epsilon^{-2}$, we have that 
\begin{equation}
    \min_{k=0,\dots, n} \e [h_k] \leq \sqrt{M \sigma^2 \e [V_0]} \cdot \epsilon,
\end{equation}
which completes the proof.

\end{proof}

\subsection{Proofs in Section \ref{sec: 4.4}}

\begin{proof}

Note that $\alpha = \frac{1}{L_2}$ and $\eta = \frac{18 p}{1 - p} \kappa^2 \alpha$.
Without loss of generality, we assume that $\inf_x V(x) = \inf V > - \infty$.
Recall the explicit of equation \eqref{equ: a26}, for any $k \geq 0$, we have
\begin{equation} \label{equ: a63}
   h_k \leq \frac{1}{\alpha p} (V_k - \e_k [V_{k+1}]) + \frac{18^2 p}{2 (1-p)} \kappa^4 L_1 \alpha \sigma^2 + \mathrm{const}, 
\end{equation}
where $\mathrm{const} = \frac{11}{20} L_2 \alpha^2 + \frac{1}{2} L_1 \alpha^2$.
Taking expectation to \eqref{equ: a63} and summing for $k = 1, 2, \dots, n$, we have 
\begin{equation} \label{equ: 65}
    \sum_{k=1}^n h_k \leq \frac{1}{\alpha p} (V_1 - \inf V) + \frac{18^2 p n}{2 (1-p)} \kappa^4 L_1 \alpha \sigma^2 + n \cdot \mathrm{const}.
\end{equation}
Hence, to minimize $\sum_{k=1}^n h_k$, it is sufficient to minimize $\frac{1}{\alpha p} (V_1 - \inf V) + \frac{18^2 p n}{2 (1-p)} \kappa^4 L_1 \alpha \sigma^2$.
Let $\delta = V_1 - \inf V$, and the RHS of \eqref{equ: 65} is denoted by $\Theta (p)$.

First, note that the assumption $\eta \leq \frac{1}{L_1}$ is equivalent to
\begin{equation}
    p \leq \frac{L_2}{9 L_1 \kappa^2 + L_2}.
\end{equation}

\underline{\textbf{Case 1}}. $\sigma = 0$.

The minimization problem is trivial in this case.
We have $\arg\min_p \Theta (p) = \frac{L_2}{9 L_1 \kappa^2 + L_2}$.

\underline{\textbf{Case 2}}. $\sigma > 0$.

On the one side, the minimization problem $\arg\min_{p > 0} \left\{ \frac{1}{\alpha p} \delta + \frac{18^2 p n}{2 (1-p)} \kappa^4 L_1 \alpha \sigma^2 \right\}$ is equal to
\begin{equation}
    p = \frac{\sqrt{\delta} \sqrt{\delta + 648 \alpha^2 \kappa^4 L_1 \sigma^2 n} - \delta}{324 \alpha^2 \kappa^4 L_1 \sigma^2 n}.
\end{equation}

On the other side, note that the assumption $\eta \leq \frac{1}{L_1}$ is equivalent to
\begin{equation}
    p \leq \frac{L_2}{9 L_1 \kappa^2 + L_2}.
\end{equation}

Let $p_1 = \frac{\sqrt{\delta} \sqrt{\delta + 648 \alpha^2 \kappa^4 L_1 \sigma^2 n} - \delta}{324 \alpha^2 \kappa^4 L_1 \sigma^2 n}$ and $p_2 = \frac{L_2}{9 L_1 \kappa^2 + L_2}$.
According to the basic theory of optimization, we have the following result.
If $p_1 \leq p_2$, we have $ \mathop{\arg\min}_{p \in (0, 1)} \Theta (p)$.
Else, we have $\mathop{\arg\min}_{p \in (0, 1)} \Theta (p) = p_2$.
Hence, combining the results above together, we have 
\begin{equation}
    \mathop{\arg\min}_{p \in (0, 1)} \Theta (p) = \min \left\{ p_1, p_2 \right\}.
\end{equation}

\end{proof}
\end{appendices}

\end{document}